\documentclass[12pt]{article}

\usepackage{graphicx,float,latexsym,times}
\usepackage{amsfonts,amstext,amsmath,amssymb,amsthm}
\usepackage{caption2}
\usepackage{color}

\long\def\symbolfootnote[#1]#2{\begingroup%
\def\thefootnote{\fnsymbol{footnote}}\footnote[#1]{#2}\endgroup}
\usepackage{titlesec} % Very fancy section title control: Sets sections as SQA requires

\titleformat{\section}{\large\bfseries}{\thesection.}{.5em}{}
\titlespacing*{\section}{0pt}{*3}{*2}
\titleformat{\subsection}{\normalfont\bfseries}{\thesubsection.}{.5em}{}
\titlespacing*{\subsection} {0pt}{*3}{*2}
\titleformat{\subsubsection}{\normalfont\bfseries}{\thesubsubsection.}{.5em}{}
\titlespacing*{\subsubsection} {0pt}{*3}{*2}

\theoremstyle{plain} %% italic text
\newtheorem{theorem}{Theorem}[section]
\newtheorem{lemma}{Lemma}[section]

\newtheorem{proposition}{Proposition}[section]
\newtheorem{example}{Example}[section]

\theoremstyle{definition} %% or \theoremstyle{remark} will produce roman text

\newtheorem{remark}{Remark}[section]

\numberwithin{equation}{section} %% double numbering within sections

\usepackage[authoryear]{natbib}
\definecolor{mred}{rgb}{0, 0, 0}
\definecolor{dgreen}{rgb}{0,0,0}
\definecolor{mgreen}{rgb}{0,0.75,0}
\definecolor{dblue}{rgb}{0,0,0}
\definecolor{mblue}{rgb}{0,0,0.75}
\definecolor{dyellow}{rgb}{0.5,0.5,0}
\definecolor{dmagenta}{rgb}{0.75,0,0.75}
\definecolor{dcyan}{rgb}{0,0.5,0.5}
\def\tmred{\textcolor{mred}}

\def\tdgreen{\textcolor{dgreen}}

\def\tdblue{\textcolor{dblue}}

\begin{document}

\title{\textbf{\Large Controlled Sensing for Sequential Multihypothesis Testing\\ with Controlled Markovian Observations and Non-Uniform Control Cost}}

\date{}

\maketitle

%%%%%%%%% Authors, affiliations %%%%%%%%%%%%%%%%%%%%%%%%%%

\author{
\begin{center}
\vskip -1cm
\textbf{\large Sirin Nitinawarat and Venupogal V. Veeravalli}\\
Department of Electrical and Computer Engineering\\
and\\
Coordinated Science Laboratory\\
University of Illinois at Urbana-Champaign\\
Urbana, IL 61801, USA
\end{center}
}

\symbolfootnote[0]{\normalsize Address correspondence to Venugopal V. Veeravalli, Department of ECE, University of Illinois at Urbana-Champaign, 106 Coordinated Science Laboratory, 1308 West Main Street, Urbana, IL 61801, USA; E-mail: vvv@illinois.edu}

{\small \noindent\textbf{Abstract:  A new model for controlled sensing for multihypothesis testing is proposed and studied in the sequential setting.  This new model, termed {\em controlled Markovian observation} model, exhibits a more complicated memory structure in the controlled observations than existing models.  In addition, instead of penalizing just the delay until the final decision time as \tdgreen{in} standard sequential hypothesis testing problems, a much more general cost structure is considered which entails accumulating the total control cost with respect to an arbitrary \tdgreen{control cost function.}  An asymptotically optimal test is proposed for this new model and is shown to satisfy 
\tdgreen{an {\em optimality}} condition formulated in terms of decision making risk.  It is shown that the optimal causal control policy for the controlled sensing problem is self-tuning, \tdgreen{in the sense of}  maximizing an inherent ``inferential'' reward simultaneously under every hypothesis, with the maximal value being the best possible \tdgreen{corresponding to the case where} the true hypothesis is known at the outset.  Another test is also proposed to meet {\em distinctly predefined} constraints on the various decision risks {\em non-asymptotically,} while retaining asymptotic optimality.
} }
\\ \\
%%%%%%%%% Key words %%%%%%%%%%%%%%%%%%%%%%%%%%
\\ \\
%%%%%%%%% Subject Classifications %%%%%%%%%

\section{INTRODUCTION} \label{s:Intro}

The broad topic of controlled sensing for inference primarily concerns adaptively managing and controlling multiple degrees of freedom (e.g., sensing modality) in an information-gathering system to solve a given inference task.  Contrary to the traditional control systems wherein the control primarily affects the evolution of the state, in controlled sensing for inference the control affects only the observations.  The goal in controlled sensing for inference is {\em not} to drive the state to the desired range, but to best shape the quality of observations to facilitate solving the given inference task.

Applications of controlled sensing for inference include, but are by no means limited to, adaptive resource management in sensor and wireless networks \citep{hero-etal-sens-book-2008, bane-veer-sqa-2012}, medicine \citep{ansc-jasta-1963}, and search \citep{ston-search-book-1989} and target tracking \citep{fuem-veer-ieeetsp-2008}.

In this paper, we focus on the basic inference problem of sequential hypothesis testing, and our goals are first, to design an asymptotically optimal test consisting of a causal control policy, a stopping rule, and a final decision rule; and second, to characterize its asymptotic performance.  Special emphasis will be placed on characterizing the asymptotically optimal control policy.

Controlled sensing for sequential multihypothesis testing has been previously studied for the model in which, conditioned on the control signal at each time, the observation at that time is conditionally independent of all past observations and control signals.  Under a suitable ``positivity'' assumption imposed on the model, an asymptotically optimal test using a {\em stationary causal control} policy was proposed, first by \cite*{cher-amstat-1959} for binary hypothesis testing, and then by \cite*{bess-tech-repo-1960} for general multihypothesis testing.  \tdgreen{In these papers, the authors termed the considered problem ``sequential design of experiments.''}
With the positivity assumption, \tdgreen{the asymptotically optimal test of Chernoff and Bessler} was shown to achieve the smallest asymptotes of expected \tdgreen{stopping times} under every hypothesis simultaneously as the probability of error vanishes. Another test was proposed in \cite*{niti-atia-veer-ieeetac-2013} that can be used to successfully  dispense with the positivity assumption imposed in \cite*{cher-amstat-1959, bess-tech-repo-1960}.  It was also proven in \cite*{niti-atia-veer-ieeetac-2013} that the proposed test satisfies \tdblue{a different asymptotic optimality condition, formulated in terms of decision making risk,} than the one adopted in \cite*{cher-amstat-1959, bess-tech-repo-1960}.  In \cite{nagh-javi-astat-2013}, a Bayesian version of the controlled sensing for sequential multihypothesis testing problem was considered in the non-asymptotic regime, and the dynamic programming equation was analyzed to find the structure of the optimal test.

Controlled sensing for sequential multihypothesis testing is related to but \tdgreen{is also} fundamentally different from feedback channel coding with variable-length codes \tdgreen{introduced in} \citep{burn-prob-per-inf-1976}.  First, the controller in the controlled sensing problem does {\em not} know the hypothesis, whereas the feedback encoder in the channel coding problem {\em knows} the message.  Second, the distributional model of the observation in the controlled sensing problem can depend on the hypothesis in an {\em arbitrary} manner, whereas in the channel coding problem, the distribution of the channel output is constrained to depend on the message through a {\em fixed} channel that is {\em independent} of the message.  These two distinctions make the controlled sensing problem more challenging than the feedback channel coding problem when the number of messages is {\em fixed,} and, hence, make the vast literature of feedback channel coding with variable-length codes (see, e.g., \cite*{burn-prob-per-inf-1976, yama-itoh-ieeetit-1979, ooi-worn-ieeetit-1998, naki-gall-ieeetit-2008, como-yuks-tati-ieeetit-2009}) not directly relevant to the controlled sensing problem.  Having said that, we note that we only consider controlled sensing models with a fixed number of hypotheses, whereas the main focus in the channel coding problem is for an exponentially growing number of messages (with a fixed rate).  Consequently, in the regime of a growing number of messages, the feedback channel coding problem is also unique from our controlled sensing problem.

In this paper, we propose a new observation model for controlled sensing for sequential multihypothesis testing termed a {\em controlled Markovian observation} model.  The memory structure in the controlled observations for this new model is more complicated than in the existing model considered in all the previous work \citep{cher-amstat-1959, bess-tech-repo-1960, niti-atia-veer-ieeetac-2013, nagh-javi-astat-2013}.  In addition, we consider a more general cost structure than just the delay until the final decision as in all the previous work.  In particular, we allow for any total accumulated control cost with respect to an arbitrary control cost function, up to the final decision time.
%MOTIVATION TO CONSIDER THE MORE GENERAL COST STRUCTURE

Our technical contributions are as follows.  We propose an asymptotically optimal test for this new observation model and show that it satisfies 
\tdblue{an {\em asymptotic optimality}} condition formulated in terms of decision making risk.  We also propose another test that meets {\em distinctly predefined} constraints on the various decision risks {\em non-asymptotically,} while retaining asymptotic optimality.

Our results are proven using a combination of tools and principles from both decision theory and stochastic control.  Interestingly, although the role of the causal control policy in the controlled sensing problem is merely to facilitate the eventual testing among the hypotheses without any explicit reward structure to gauge how well the different control policies perform, our results show that there is an {\em inherent} ``inferential'' reward structure maximized by the control policy of the asymptotically \tdblue{optimal} test for the controlled sensing problem.  This reward is given in terms of \tdgreen{the}  ratio between the time average of a suitable decision-theoretic distance of the controlled observations under the true hypothesis, to all competing hypotheses, and the time average of the control cost.  More importantly, we show that the optimal causal control policy for the controlled sensing problem achieves the optimal limiting value of this inferential reward simultaneously under {\em every} hypothesis.  In addition, this optimal limiting value is also the best possible in the fictitious situation where the true hypothesis is known to the controller at the outset.  In effect, the optimal causal control policy is able to {\em self-tune} to the true hypothesis while achieving the optimal limit of the inferential reward.  Such a control policy is reminiscent of the {\em self-tuning} optimal control policy extensively studied in adaptive stochastic control (see, e.g., \cite*{kuma-vara-ss-book-1986}).

\tdblue{The material in this paper was presented in parts at the IEEE International Symposium on Information Theory, Istanbul, Turkey, July 7-12, 2013, and also at the 2013 Asilomar Conference on Signals, Systems and Computers, November 3-6, 2013.}

\section{PRELIMINARIES AND MODELS WITH UNCONTROLLED\\ OBSERVATIONS}
\label{s:Prelim&UncontrolledObservations}

Throughout this paper, all random variables are assumed to take values in \tdblue{{\em finite} sets} and are denoted by capital letters; their realizations are denoted by the corresponding lower-case letters.

We start by considering sequential multihypothesis testing based on observation sequence which cannot be controlled.  Let the joint distribution of the observations $Y_1, \ldots, Y_k,\ k \geq 1,$ under hypothesis $\,i = 0, \ldots, M-1\,$ be denoted by \tdblue{a pmf} $\,p_i \left( y^k \right),$ where $\,y^k\,$ denotes $\,\left( y_1, \ldots, y_k \right).$

A sequential test consists of a stopping rule and a final decision rule.  The stopping rule determines a (random) number of observations, denoted by $N$, that are taken until a decision is made for a hypothesis according to $\,\delta \left( Y^N \right)$.  The overall goal of a sequential test is to make the final decision with a desirable level of accuracy using the fewest number of observations.

There are a number of metrics that can be used to gauge the accuracy of a sequential test.  The most conservative one is the worst-case (maximal) probability of error $\,P_{ {\scriptsize \mathrm{max}}}\,$ defined as
\begin{equation}
	P_{ {\scriptsize \mathrm{max}}} \ \triangleq\ 
	\max_{i=0, \ldots, M-1}~\mathbb{P}_i \left \{ \delta \neq i \right \}
	\ =\ 
	\max_{i=0, \ldots, M-1}~\mathbb{P}_i \left \{ \delta \left( Y^N \right) \neq i \right \},
	\label{eqn-def-Pmax}
\end{equation}
where each $\,\mathbb{P}_i \left \{ \delta \neq i \right \}\,$ is the conditional probability of error under hypothesis $\,i = 0, \ldots, M-1$.  Another more refined notion is that of \tdblue{maximal decision risks} each of which $\,R_i,\ i = 0, \ldots, M-1,\,$ is defined to be the \tdblue{maximal conditional probability} of incorrectly deciding for hypothesis $\,i\,$ as
\tdblue{
\begin{equation}
	R_i \ \triangleq\ 
	\max\limits_{j \neq i}~\mathbb{P}_{j} \left \{ \delta = i \right \}.
	\label{eqn-def-Risks}
\end{equation}
}
%where $\,\pi\,$ is the prior distribution of the hypothesis.  
It clearly follows from (\ref{eqn-def-Pmax}) and (\ref{eqn-def-Risks}) that
for every $\,i = 0, \ldots, M-1,$
\begin{equation}
	R_i \ \leq\ P_{ {\scriptsize \mathrm{max}}}.
	\label{eqn-Pmax-dom-Risks}
\end{equation}
\tdblue{We denote the class of tests that satisfy $P_{ {\scriptsize \mathrm{max}}} \leq \alpha;$ by $\mathbb{C} \left( \alpha \right);$ and that satisfy $R_0 \leq \beta_0, \ldots, R_{M-1} \leq \beta_{M-1}$ by $\mathbb{C}' \left( \beta_0, \ldots, \beta_{M-1} \right),$ respectively.}
In the following presentations, our results will be stated in terms of the Kullback-Leibler (KL) distance between two distributions (pmfs) $\,p_1(y), p_2(y)\,$ on a common finite \tdblue{set} $\,\mathcal{Y},$ which is defined as
\begin{equation}
	D \left( p_1 \| p_2 \right)
	\ \triangleq\ \sum\limits_{y \in \mathcal{Y}}
	~p_1(y) \log{\left( \frac{p_1(y)}{p_2(y)} \right)},
	\label{eqn-KLdistance}
\end{equation}
if the support of $\,p_2\,$ subsumes that of $\,p_1$ 
\tdblue{(with the convention that $0 \log{0} = 0$).}
%; and is defined to be infinite otherwise.

\subsection{Asymptotically Optimal Test: Multihypothesis Sequential Probability Ratio Test (MSPRT)}
\label{ss:MSPRT}

We now describe a sequential test that is asymptotically optimal in the regime wherein the maximal probability of error vanishes (hence by (\ref{eqn-Pmax-dom-Risks}), so do all the risks).  The stopping rule of this test entails stopping at the first time at which the likelihood ratios between the maximum likelihood estimate of the hypothesis and all other competing hypotheses are larger than a threshold $\,T > 0,$ i.e.,
\begin{equation}
\tdblue{
	N^* \ \triangleq\ \min \left \{n \geq 1 \ \Bigg\vert\ 
	\frac{p_{\hat{i}} \left(Y^n \right)}{\max\limits_{j \neq \hat{i}}~p_j \left( Y^n \right)}
	~>~T \right \},
\label{eqn-MSPRT}
}
\end{equation}
where $\,\hat{i} ~\triangleq~ \mathop{\mbox{argmax}}\limits_{i = 0, \ldots, M-1}~p_i \left( Y^n \right).$  At the stopping time a decision is made for the most likely hypothesis: \tdblue{$\delta^* \left( Y^{N^*} \right) = \hat{i}.$}  
\tdblue{We denote the error probability and risks for the test $\left( N^*, {\delta}^* \right)$ by
$P_{ {\scriptsize \mathrm{max}}}^*$ and $R_0^*, \ldots R_{M-1}^*,$ respectively.}
The stopping rule (\ref{eqn-MSPRT}) was first proposed by Chernoff \citep{cher-amstat-1959} in the setup with observation control.  A similar test was later proposed for models with uncontrolled observations \citep{baum-veer-ieeetit-1994} and was termed ``Multihypothesis Sequential Probability Ratio Test'' (MSPRT) therein.  In this other test, all the likelihoods in (\ref{eqn-MSPRT}) are weighted by the respective prior probabilities, and the maximum in the denominator is changed to summation.  Another similar test was also proposed in \cite*{drag-et-al-ieeetit-1999} with the weighting by the prior probabilities while retaining the maximum in the denominator.  Since all of these tests are asymptotically optimal when the observations are independent and identically distributed (i.i.d.) but the test in (\ref{eqn-MSPRT}) does not require the knowledge of the prior distribution 
%$\,\pi,$ 
\tdblue{of the hypothesis.}
we choose to focus on the test in (\ref{eqn-MSPRT}) and call it MSPRT.

\begin{proposition}
[\tdblue{Asymptotic Optimality of the MSPRT \citep{veer-baum-ieeetit-1995}}]
\label{prop1-uncontrolled}
Assume that the observations under each hypothesis are i.i.d. according to $\,p_i(y),\ i = 0, \ldots, M-1,$ and that it holds for every \tdblue{$\, i, j = 0, \ldots, M-1,\ i \neq j,$ that $\,D \left( p_i \| p_j \right) ~<~ \infty.$}  \tdblue{Then any sequence of tests with vanishing maximal probability of error, $P_{ {\scriptsize \mathrm{max}}} \ \rightarrow\ 0,$ satisfies for every $\,i = 0, \ldots, M-1$ that 
\begin{equation}
\mathbb{E}_i [N] \ \geq\ 
\frac{-\log{\left( 
R_i
\right)}}
{\min\limits_{j \neq i}~D \left( p_i \| p_j \right)} \left( 1+ o(1) \right),
\label{eqn-converse-uncontrolled1}	
\end{equation}
hence, also that
\begin{equation}
\mathbb{E}_i [N] \ \geq\ 
\frac{-\log{\left( 
P_{ {\scriptsize \mathrm{max}}}
\right)}}
{\min\limits_{j \neq i}~D \left( p_i \| p_j \right)} \left( 1+ o(1) \right).
\label{eqn-converse-uncontrolled2}
\end{equation}}

\tdblue{For the MSPRT, it holds that 
\begin{equation}
\tdblue{P_{ {\scriptsize \mathrm{max}}}^* \ \leq \frac{M-1}{T}}.
\end{equation}  
Furthermore, the MSPRT is \tdblue{asymptotically optimal.}  
In particular, as the threshold $\,T\,$ in (\ref{eqn-MSPRT}) approaches infinity, the MSPRT yields for each $\,i = 0, \ldots, M-1,\,$ that
\begin{equation}
\mathbb{E}_i [N^*]\ =\ 
\frac{-\log{\left( 
P_{ {\scriptsize \mathrm{max}}}^* 
\right)}}
{\min\limits_{j \neq i}~D \left( p_i \| p_j \right)} \left( 1+ o(1) \right)
\ =\ 
\left( \inf_{\left(N, \delta \right) \,\in\, \mathbb{C} \left( P_{{\scriptsize \mathrm{max}}}^* \right)} 
\mathbb{E}_i [N] \right)
(1+o(1)),
\label{eqn-perf-MSPRT1}
\end{equation}
hence, also that
\begin{equation}
\mathbb{E}_i [N^*]\ =\ 
\frac{-\log{\left( 
R_i^*
\right)}}
{\min\limits_{j \neq i}~D \left( p_i \| p_j \right)} \left( 1+ o(1) \right)
\ =\ 
\left( 
\inf_{\left(N, \delta \right) \,\in\, \mathbb{C}' \left( R_0^*, \ldots, R_{M-1}^* \right)} \mathbb{E}_i [N] \right)
(1+o(1)).
\label{eqn-perf-MSPRT2}
\end{equation}}

\end{proposition}
\begin{remark}
In the forward claim above, it suffices just to establish \tdblue{the first equality in} (\ref{eqn-perf-MSPRT1}), as \tdblue{the first equality in} (\ref{eqn-perf-MSPRT2}) follows immediately from (\ref{eqn-perf-MSPRT1}) using (\ref{eqn-Pmax-dom-Risks}).  \tdblue{The second equalities in (\ref{eqn-perf-MSPRT1}) and (\ref{eqn-perf-MSPRT2}) follow from the first ones and (\ref{eqn-converse-uncontrolled1}), (\ref{eqn-converse-uncontrolled2}), respectively.}  
In the reverse direction, only the assertion (\ref{eqn-converse-uncontrolled1}) needs to be established.  That the asymptotic performances with respect to both the maximal probability of error and the corresponding risks coincide (the optimal coefficients in (\ref{eqn-perf-MSPRT1}), (\ref{eqn-converse-uncontrolled2}) and in (\ref{eqn-perf-MSPRT2}), (\ref{eqn-converse-uncontrolled1}) are the same) is what we \tdblue{mean by {\em asymptotic optimality.}}
\end{remark}

It is interesting to note that the aforementioned result pertaining to the \tdblue{asymptotic optimality} relies critically on the notion of decision risks.  \tdblue{In particular, we now show in the following simple counterexample that the version of this claim with respect to the more natural criteria of conditional probabilities of error {\em does not hold}.}  Note that each of the conditional probabilities of error is still dominated by the maximal probability of error similar to each of the risks as in (\ref{eqn-Pmax-dom-Risks}).  The characterization of the optimal asymptotic performance with respect to the conditional probabilities of error is still unresolved.

\tdblue{
\begin{example}  Consider binary hypothesis testing between two distributions $p_0, p_1$ on a finite $\mathcal{Y}$ with 
$D \left( p_0 \| p_1 \right) < \infty,\ D \left( p_1 \| p_0 \right) < \infty.$  For a fixed constant $c>0$ and a threshold $T> 1,$ consider the test with the stopping rule being
\begin{align}
	N = \min \left(
		N_0 = \mathop{\rm{argmin}}_{n \geq 1} 
			\left( \frac{p_0 \left( Y^n \right)}{p_1 \left( Y^n \right)} \geq T^c \right),
		N_1 = \mathop{\rm{argmin}}_{n' \geq 1},
			\left( \frac{p_1 \left( Y^{n'} \right)}{p_0 \left( Y^{n'} \right)} \geq T \right)
	\right);
	\nonumber
\end{align}
and with the ML rule as the final decision rule.
Then, it can be shown that as $T \rightarrow \infty,$
\begin{align}
	\frac{N_0}{\log{T}}   &\ \stackrel{\rm{a.s.}}{\rightarrow}\ 	
		\frac{c}{D \left( p_0 \| p_1 \right)}\ \ \ \ \mbox{under}~\mathbb{P}_0;\nonumber \\
	\frac{N_1}{\log{T}}   &\ \stackrel{\rm{a.s.}}{\rightarrow}\ 	
		\frac{1}{D \left( p_1 \| p_0 \right)}\ \ \ \ \mbox{under}~\mathbb{P}_1;\nonumber \\
	\mathbb{P}_0 &\ \left \{ \delta = 1 \right \}  \leq\ \frac{1}{T};	\nonumber \\
	\mathbb{P}_1 &\ \left \{ \delta = 0 \right \}  \leq\ \frac{1}{T^c}.  \nonumber
\end{align}
Consequently, we get that as $T \rightarrow \infty,$
\begin{align}
	\mathbb{E}_0 \left [N \right] 
	&\ \leq\ \frac{-\log{\mathbb{P}_0 \left \{ \delta = 1 \right \}}}{\frac{D \left( p_0 \| p_1 \right)}{c}}
			(1 + o(1))
	\label{eqn-ex1-3} \\
	\mathbb{E}_1 \left [N \right] 
	&\ \leq\ \frac{-\log{\mathbb{P}_1 \left \{ \delta = 0 \right \}}}{ D \left( p_1 \| p_0 \right) c}
			(1 + o(1)).
	\label{eqn-ex1-4}
\end{align}
By selecting $c > 1,$ we can trade the asymptotic performance under hypothesis $0$ (characterized by the number $\frac{D \left( p_0 \| p_1 \right)}{c} < D \left( p_0 \| p_1 \right)$) off that under hypothesis $1$ (characterized by the number $c D \left( p_1 \| p_0 \right) > D \left( p_1 \| p_0 \right)$), and vice versa, by selecting $c, 0< c< 1$.  In contrary, if we consider the asymptotic performance in terms of decision risks, then it follows from Proposition \ref{prop1-uncontrolled} that for any sequence of tests it must hold simultaneously that
\begin{align}
	\mathbb{E}_0 \left [N \right] 
	&\ \geq\ \frac{-\log{R_0}}{D \left( p_0 \| p_1 \right)} (1 + o(1))
	\ =\ \frac{-\log{\mathbb{P}_1 \left \{ \delta = 0 \right \}}}{D \left( p_0 \| p_1 \right)} (1 + o(1))
	\label{eqn-ex1-1} \\
	\mathbb{E}_1 \left [N \right] 
	&\ \geq\ \frac{-\log{R_1}}{D \left( p_1 \| p_0 \right)} (1 + o(1))
	\ =\ \frac{-\log{\mathbb{P}_0 \left \{ \delta = 1 \right \}}}{ D \left( p_1 \| p_0 \right)} (1 + o(1)).
	\label{eqn-ex1-2}
\end{align}
In other words, there cannot exist a sequence of tests with the asymptotic performance beating either (\ref{eqn-ex1-1}) (with the corresponding number in the denominator being larger than 
$D\left( p_0 \| p_1 \right)$) or (\ref{eqn-ex1-2}) (with the corresponding number in the denominator being larger than $D\left( p_1 \| p_0 \right)$).\\
\indent
%Although the tradeoff between the two relationships in (\ref{eqn-ex1-3}), (\ref{eqn-ex1-4}) between conditional probabilities of error and the expected values of the stopping time under the corresponding hypotheses should be asymptotically optimal, the situation becomes much more difficult to resolve with more than two hypotheses.
\end{example}
}

%\begin{example}  
%\label{ex-1}
%Consider an example with three hypotheses and a ternary $\,\mathcal{Y} = \left \{ 0, 1, 2 \right \}.$  In particular, let $\,p_0(0) = 0, p_0(1) = p_0 ( 2 )= \frac{1}{2};\  p_1 (0) = p_1 (1) = \frac{1}{4}, p_1 ( 2 ) = \frac{1}{2};$ and $\,p_2 (0) = p_2 ( 2 ) = \frac{1}{4},\ p_2 ( 1 ) = \frac{1}{2}.$  It is clear that 
%\begin{equation}
%	\min \left( D \left( p_0 \| p_1 \right), D \left( p_0 \| p_2 \right) \right)
%	\ =\ D \left( p_0 \| p_1 \right) \ <\ \infty.
%	\label{eqn-ex1}
%\end{equation}
%Now consider any stopping time $\,N,$ and a decision rule which makes a decision for hypothesis $\,0$ only if the observation symbol $0$ never appears.  It is then clear that $\,\mathbb{P}_0 \left \{ \delta \neq 0 \right \} = 0.$ Hence $\,\mathbb{E}_0 [N]\,$ can be made arbitrarily small relative to the conditional probability of error $\,\mathbb{P}_0 \left \{ \delta \neq 0 \right \} = 0,$ whereas the associated minimum KL distance is finite, thereby exceeding the optimal asymptotic performance with respect to the maximal probability of error and the risks in (\ref{eqn-converse-uncontrolled1}), (\ref{eqn-converse-uncontrolled2}).
%\end{example}

\section{PROBLEM FORMULATION}
\label{s:Formulation}

We now move on to the main subject of the paper, namely sequential multihypothesis testing with observation control.  The observation and control signal at each time step are assumed to take values in {\em finite} sets $\,\mathcal{Y}\,$ and $\,\mathcal{U},$ respectively.  At each time $\,k = 1, 2, \ldots,$ conditioned on each hypothesis $\,i = 0, \ldots, M-1,$ on the current control value $\,u_k\,$ and on the previous observation $\,y_{k-1},$ the current observation $\,Y_k\,$ is assumed to be conditionally independent of all earlier observations and control signals $\,\left( y^{k-2}, u^{k-1} \right) ~\triangleq~ \left( \left( y_1, \ldots, y_{k-2} \right),  \left( u_1, \ldots, u_{k-1} \right) \right)$.  In addition, conditioned on the true hypothesis $\,i\,$ and on $\, \left( u_k, y_{k-1} \right),\ Y_k\,$ is assumed to be conditionally distributed according to $\,p_i^{u_k} \left( \cdot \vert y_{k-1} \right),$ where $\left \{ p_i^u \left( \cdot \vert \cdot \right) \right \}_{i = 0, \ldots, M-1}^{u \in \mathcal{U}}\,$ is a fixed collection of transitions probabilities from $\,\mathcal{Y}\,$ to itself.  We assume throughout that for every $\,i = 0, \ldots, M-1,\ u \in \mathcal{U},\ y, \tilde{y} \in \mathcal{Y},$
\begin{equation}
	p_i^u \left( y \vert \tilde{y} \right) \ >\ 0.
	\label{eqn-assump-pos-transition}
\end{equation}
The initial observation $\,y_0\,$ is assumed to be a constant.

Controlled sensing for sequential hypothesis testing has been previously studied for a simpler observation model \citep{cher-amstat-1959, bess-tech-repo-1960, niti-atia-veer-ieeetac-2013, nagh-javi-astat-2013}  in which for every $\,u \in \mathcal{U},\ i = 0, \ldots, M-1,$ the $\,p_i^u \left( \cdot \vert \tilde{y} \right)\,$ in (\ref{eqn-assump-pos-transition}) is independent of $\,\tilde{y}.$  Hence, for this simpler model, conditioned on each hypothesis $\,i\,$ and $\,u_k,\ Y_k\,$ is conditionally independent of $\,\left( y^{k-1}, u^{k-1} \right)$.

Our main interest here will be on {\em causal} control policies for which the control $\,U_k\,$
at each time $\,k = 1, 2, \ldots,$ can be any (possibly randomized) function of past observations and past control signals.  In particular, each $\,U_k,\ k = 2, 3, \ldots,$ is specified by an arbitrary conditional pmf $\,q_k \left( u_k \vert y^{k-1}, u^{k-1} \right),$ and $\,U_1\,$ is specified by a conditional pmf $\,q_1 \left( u_1 \vert y_0 \right).$  Having specified the causal control policy, the joint distribution of observations and control signals under hypothesis $\,i,$ denoted by $\,p_i \left( y^n, u^n \right),$ can be written as 
\begin{equation}
	p_i \left( y^n, u^n \right)
	\ \triangleq\ 
	q_1 \left( u_1 \vert y_0 \right) 
	\prod_{k=1}^n  p_i^{u_k} \left( y_k \vert y_{k-1} \right)
	\prod_{k=2}^n  q_k \left( u_k \vert y^{k-1}, u^{k-1} \right). 
	\label{eqn-jointdist-obs-control}
\end{equation}
The entire collection of conditional pmfs $\,\left \{ q_1 \left( u_1 \vert y_0 \right), 
\left \{ q_k \left( u_k \vert y^{k-1}, u^{k-1} \right) \right \}_{k=2}^{\infty} \right \}\,$ describes a causal control policy which will be denoted collectively by $\,\phi$.

Let $\,\mathcal{F}_k\,$ denote the $\,\sigma$-field generated by $\, \left( Y^k, U^k \right).$  A {\em sequential test} $\, \gamma = \left( \phi, N, \delta \right)\,$ consists of a causal control policy $\,\phi,$ an $\,\mathcal{F}_k$-stopping time $\,N\,$ denoting the (random) number of observations up to the time of the final decision, and a decision rule $\,\delta \left( Y^N, U^N \right)\,$ selecting one of the hypotheses in $\,\left \{ 0, \ldots, M-1 \right \}$.

\section{ASYMPTOTICALLY OPTIMAL TEST AND ITS PERFORMANCE}
\label{s:AsymptoticallyOptimalTest}

\subsection{\tdblue{Proposed Test} and Its \tdblue{Asymptotic Optimality}}
\label{ss:TestNoCost}

%[Motivation of our test.  We stick with the design principle of using the principle underlying %MSPRT on the induced joint distribution of observations and controls to get the stopping rule.]
For every $\,i = 0, \ldots, M-1,$ and every conditional distribution $\,q \left( u \vert \tilde{y} \right),$ we define $\,p_i^{q}\,$ as the following transition probabilities from $\,\mathcal{Y}\,$ to itself
\begin{equation}
	p_i^{q} \left( y \vert \tilde{y} \right) 
	\ \triangleq\ 
	\sum\limits_{u \in \mathcal{U}}  q \left( u \vert \tilde{y} \right) p_i^u \left( y \vert \tilde{y} \right).
\end{equation}
It follows from Assumption (\ref{eqn-assump-pos-transition}) and finiteness of $\,\mathcal{Y}, \mathcal{U}\,$ that, for every $\,i, q,$ such $\, p_i^q\,$ admits a unique stationary distribution $\,\mu_i^q \left( \tilde{y} \right)$.

The causal control policy that we shall adopt admits the following sequential description.  Except for suitably chosen sparse and recurrent occasions to be specified shortly, the control signal at time $\,k\,$ is selected to be a (randomized) function of $\,y_{k-1},$ and the ML estimate of the hypothesis $\,\hat{i}_{k-1} = \mathop{\mbox{argmax}}\limits_{i=0, \ldots, M-1} p_i \left( y^{k-1}, u^{k-1} \right),$ according to
\begin{equation}
	U_k \ \sim\ 
	q_{\hat{i}_{k-1}}^* \left( \cdot \vert y_{k-1} \right),
	\label{eqn-opt-control-exploit}
\end{equation} 
where for each $\,i = 0, \ldots, M-1,$
\begin{equation}
q_i^* \left( \cdot \vert \cdot \right) 
\ =\ \mathop{\mbox{argmax}}_{q \left( \cdot \vert \cdot \right)}
~\min\limits_{j \neq i}
~\sum\limits_{\tilde{y} \in \mathcal{Y}, u \in \mathcal{U}}
\mu_i^q \left( \tilde{y} \right) q \left( u \vert \tilde{y} \right)
\tdblue{
D \left( p_i^u \left( \cdot \vert \tilde{y} \right) \| p_j^u \left( \cdot \vert \tilde{y} \right) \right).
}
\label{eqn-opt-exploit-control-nocost}
\end{equation}
The control (\ref{eqn-opt-control-exploit}) is followed at all times except at time steps
\begin{equation}
k= \lceil a^{\ell} \rceil,\ \ell = 0, 1, \ldots, 
\label{eqn-sparse-exploration}
\end{equation}
when $\,U_k\,$ is picked to be uniformly distributed independently of the the ML estimate of the hypothesis and the entire past observations and past control signals.  In (\ref{eqn-sparse-exploration}), $\,a > 1$ is picked to be sufficiently close to 1 to make these occasions in (\ref{eqn-sparse-exploration}) not too sparse.  \tdblue{Heuristically speaking, the control rule in (\ref{eqn-sparse-exploration}) has to {\em explore} often enough (not too sparse) in order to ensure that the ML estimate of the hypothesis converges to the true hypothesis sufficiently quickly.  On the other hand, the other control rule in (\ref{eqn-opt-control-exploit}), applied for most of the times, exploits the ML estimate of the hypothesis to select the most informative observation in order to yield the asymptotically optimal performance.}

As for the stopping rule, we still apply the same principle as in the MSPRT except that we now instead use the joint distributions under various hypotheses of all the observations and control signals up to the current time; these joint distributions are induced by the causal control policy (\ref{eqn-opt-control-exploit}), (\ref{eqn-opt-exploit-control-nocost}), (\ref{eqn-sparse-exploration}) according to (\ref{eqn-jointdist-obs-control}).  In particular, we stop at the first time $\,N\,$ at which
\begin{equation}
\tdblue{
N^* \ \triangleq\ \min \left \{ n \geq 1  \  \Bigg \vert\ 
	\frac{p_{\hat{i}_n} \left(Y^n, U^n \right)}{\max\limits_{j \neq \hat{i}_n}~p_j \left( Y^n, U^n \right)}
	~>~T \right \}.
	\label{eqn-stopping-nocost}
}
\end{equation}
At the stopping time a decision is made for the most likely hypothesis \tdblue{with the incurred error probability and risks being $P_{ {\scriptsize \mathrm{max}}}^*$ and $R_0^*, \ldots R_{M-1}^*,$ respectively.}

\begin{theorem}[\tdblue{Asymptotic Optimality}]  
\label{thm-optasymperf-nocost}
\tdblue{
Any sequence of tests with vanishing maximal probability of error, $P_{ {\scriptsize \mathrm{max}}} \ \rightarrow\ 0,$ satisfies for every $\,i = 0, \ldots, M-1\,$ that 
\begin{eqnarray}
\mathbb{E}_i [N] &\geq& 
\frac{-\log{\left( 
R_i
\right)}}
{
\max\limits_{q \left( \cdot \vert \cdot \right)}
~\min\limits_{j \neq i}
~\sum\limits_{\tilde{y} \in \mathcal{Y}, u \in \mathcal{U}}
\mu_i^q \left( \tilde{y} \right) q \left( u \vert \tilde{y} \right)
D \left( p_i^u \left( \cdot \| \tilde{y} \right) \| p_j^u \left( \cdot \| \tilde{y} \right) \right)
} 
\left( 1+ o(1) \right)
\label{eqn-converse-nocost1}		\\
&\geq&
\frac{-\log{\left( 
P_{ {\scriptsize \mathrm{max}}}
\right)}}
{
\max\limits_{q \left( \cdot \vert \cdot \right)}
~\min\limits_{j \neq i}
~\sum\limits_{\tilde{y} \in \mathcal{Y}, u \in \mathcal{U}}
\mu_i^q \left( \tilde{y} \right) q \left( u \vert \tilde{y} \right)
D \left( p_i^u \left( \cdot \| \tilde{y} \right) \| p_j^u \left( \cdot \| \tilde{y} \right) \right)
} 
\left( 1+ o(1) \right).
\label{eqn-converse-nocost2}
\end{eqnarray}
}

\tdblue{
For the test in (\ref{eqn-opt-control-exploit}), (\ref{eqn-opt-exploit-control-nocost}), (\ref{eqn-sparse-exploration}), (\ref{eqn-stopping-nocost}), it holds that
\begin{equation}
P_{ {\scriptsize \mathrm{max}}}^*
\ \leq \ \frac{M-1}{T}.
\end{equation}
Furthermore, the test is asymptotically optimal.  In particular, as the threshold $\,T\,$ in (\ref{eqn-stopping-nocost}) approaches infinity, \tdblue{the test} in
(\ref{eqn-opt-control-exploit}), (\ref{eqn-opt-exploit-control-nocost}), (\ref{eqn-sparse-exploration}), (\ref{eqn-stopping-nocost}) 
yields for each $\,i = 0, \ldots, M-1\,$ that
\begin{eqnarray}
\tdblue{\mathbb{E}_i [N^*]} 
&=&
\frac{-\log{\left( 
P_{ {\scriptsize \mathrm{max}}}^*
\right)}}
{
\max\limits_{q \left( \cdot \vert \cdot \right)}
~\min\limits_{j \neq i}
~\sum\limits_{\tilde{y} \in \mathcal{Y}, u \in \mathcal{U}}
\mu_i^q \left( \tilde{y} \right) q \left( u \vert \tilde{y} \right)
D \left( p_i^u \left( \cdot \| \tilde{y} \right) \| p_j^u \left( \cdot \| \tilde{y} \right) \right)
}
\left( 1+ o(1) \right)
\label{eqn-achievable-perf1}	\\
&=& 
\left( \inf_{\left( N, \delta \right) \,\in\, \mathbb{C} \left( P_{ {\scriptsize \mathrm{max}}}^* \right)} \mathbb{E}_i [N] \right)
(1+o(1))
 \nonumber \\
&=& 
\frac{-\log{\left( 
R_i^*
\right)}}
{
\max\limits_{q \left( \cdot \vert \cdot \right)}
~\min\limits_{j \neq i}
~\sum\limits_{\tilde{y} \in \mathcal{Y}, u \in \mathcal{U}}
\mu_i^q \left( \tilde{y} \right) q \left( u \vert \tilde{y} \right)
D \left( p_i^u \left( \cdot \| \tilde{y} \right) \| p_j^u \left( \cdot \| \tilde{y} \right) \right)
}
\left( 1+ o(1) \right).
\label{eqn-achievable-perf2} \\
&=& 
\left( \inf_{\left( N, \delta \right) \,\in\, \mathbb{C}' \left( R_0^*, \ldots, R_{M-1}^* \right)} \mathbb{E}_i [N] \right)
(1+o(1))
 \nonumber 
\end{eqnarray}
}
\end{theorem}

\begin{remark}  
\label{remark-nocost}
Note that the stopping time (\ref{eqn-stopping-nocost}) can be written equivalently as 
\begin{equation}
\tdblue{
N^* \ \triangleq\ 
	\min \left \{
	n \geq 1 \ \Bigg \vert\ 
	~e^{n 
	\left(
	~\min\limits_{j \neq \hat{i}_n}~\frac{1}{n}
	\log{\left( 
		\frac{p_{\hat{i}_n} \left(Y^n, U^n \right)}{p_j \left( Y^n, U^n \right)}
	\right) }
	\right )
	}
	~>~T \right \}.
	\label{eqn-stopping-nocost-alt}
}
\end{equation}
Any reasonable test would make $\,\hat{i}\,$ converge to the true hypothesis $\,i\,$ with the probability quickly converging to one.  Consequently, in order to achieve the fewest observations for the final decision-making (reaching the threshold $\,T\,$ in (\ref{eqn-stopping-nocost-alt}) as rapidly as possible), it is tempting to select a causal control policy to maximize the limit of
\begin{equation}
\min\limits_{j \neq i}~\frac{1}{n}~
	\log{\left( 
		\frac{p_{i} \left(Y^n, U^n \right)}{p_j \left( Y^n, U^n \right)}
	\right) }
\ =\ 
\min\limits_{j \neq i}~\frac{1}{n}~
\sum\limits_{k=1}^n 
	\log{\left( 
		\frac{p_{i}^{U_k} \left(Y_k \vert Y_{k-1} \right)}
		{p_{j}^{U_k} \left(Y_k \vert Y_{k-1} \right)}
	\right) }.
\label{eqn-inferential-rewards}
\end{equation}
In the fictitious situation when the true hypothesis $\,i\,$ is known to the controller at the outset, it is not hard to show using tools from Markov Decision Processes \citep{kuma-vara-ss-book-1986} that  the largest possible limiting value of (\ref{eqn-inferential-rewards}) optimized over all causal control policies (in expectation and in probability) is equal to
\begin{equation}
\max\limits_{q \left( \cdot \vert \cdot \right)}
~\min\limits_{j \neq i}
~\sum\limits_{\tilde{y} \in \mathcal{Y}, u \in \mathcal{U}}
\mu_i^q \left( \tilde{y} \right) q \left( u \vert \tilde{y} \right)
D \left( p_i^u \left( \cdot \vert \tilde{y} \right) \| 
p_j^u \left( \cdot \vert \tilde{y} \right) \right),
\label{eqn-opt-inferential-rewards}
\end{equation}
which is the coefficient in the optimal asymptotic performance in \tdblue{(\ref{eqn-converse-nocost1}), (\ref{eqn-converse-nocost2}), (\ref{eqn-achievable-perf1}), (\ref{eqn-achievable-perf2}).}  The main challenge in controlled sensing for multihypothesis testing is that the control policy has to be selected independently of the true hypothesis, as it is not known to the tester.  Despite this challenge, \tdblue{the control policy} in (\ref{eqn-opt-control-exploit}), (\ref{eqn-opt-exploit-control-nocost}), (\ref{eqn-sparse-exploration}), when perpetually applied, is able to achieve the optimal limiting value in (\ref{eqn-opt-inferential-rewards}) {\em simultaneously for every true hypothesis.}  In effect, the optimal control policy is able to {\em self-tune} to the true hypothesis while achieving the optimal limit of the ``inferential reward'' in (\ref{eqn-inferential-rewards}).  Such a control policy is reminiscent of the {\em self-tuning} optimal control policy extensively studied in adaptive stochastic control \citep{kuma-vara-ss-book-1986}.  This is also another explanation of why these coefficients in (\ref{eqn-opt-inferential-rewards}) emerge naturally in the characterization of the optimal asymptotic performance for controlled sensing for multihypothesis testing in Theorem \ref{thm-optasymperf-nocost}.  {\em Interestingly, although in the formulation of controlled sensing for multihypothesis testing there is no reward structure and the role of the control policy is merely to facilitate the eventual testing among the hypotheses, the optimal causal control policy inherently maximizes the inferential reward.}
\end{remark}

\subsection{Asymptotically Optimal Test Meeting Predefined Risk Constraints}
\label{ss:TestMeetingHardRiskConstriants}

%Although the calculation of decision risks involves, by definition (cf. (\ref{eqn-def-Risks})), the prior distribution of the hypothesis, the test proposed in Section \ref{ss:TestNoCost} does not use the knowledge of the prior distribution at all.  
\tmred{In this subsection,} we present another test that incorporates this knowledge and different thresholds for distinct identities of the ML estimate of the hypothesis.  The advantage of this new test is that it \tdgreen{enables meeting} predefined and distinct constraints on the various risks non-asymptotically \tdblue{(with potentially lower thresholds than the single threshold in (\ref{eqn-stopping-nocost}))}.  In the asymptotic regime in which all the risks vanish, we show that this new test is also \tdblue{asymptotically optimal.}  The idea of using different thresholds to meet distinctly predefined risk constraints was proposed for models with uncontrolled observations in \cite*{baum-veer-ieeetit-1994}.

Specifically, for a given tuple $\, \left( \bar{R}_0, \ldots, \bar{R}_{M-1} \right),$ we shall design a new test to satisfy $\,R_i \leq \bar{R}_i,\ i = 0, \ldots, M-1.$  To this end, we just need to modify the stopping rule (\ref{eqn-stopping-nocost}) to be 
\begin{equation}
\tdblue{
N^* \ \triangleq\ \min \left \{ n \geq 1 \ \Bigg \vert\ 
	\frac{p_{\hat{i}_n} \left(Y^n, U^n \right)}
	{\max\limits_{j \neq \hat{i}_n}~p_j \left( Y^n, U^n \right)}
	~>~ \left( \frac{1}{\bar{R}_{\hat{i}_n}} \right) \right \}.
\label{eqn-stopping-hardriskconstraints}
}
\end{equation}

\begin{theorem}[Test Meeting Distinctly Predefined Risk Constraints]  
\label{thm-optasymperf-nocost-hardrisk}
\tdblue{For any positive tuple\\ $\left( \bar{R}_0, \ldots, \bar{R}_{M-1} \right),$} %and any $\,\pi\,$ with full support, 
\tdblue{the modified test} using the stopping rule (\ref{eqn-stopping-hardriskconstraints}) instead of (\ref{eqn-stopping-nocost}) satisfies for every $\,i = 0, \ldots, M-1,$
\begin{equation}
	R_i \ \leq\ \bar{R}_i.
	\label{eqn-hard-risk-constraints}
\end{equation}
Furthermore, as $\,\max\limits_{i = 0, \ldots, M-1} \bar{R}_i \rightarrow 0,$ while satisfying $\,\max\limits_i \bar{R}_i \leq K \min\limits_i \bar{R}_i\,$ for some constant $\,K,$ this modified test is \tdblue{asymptotically optimal,} i.e., it satisfies (\ref{eqn-achievable-perf1}) and (\ref{eqn-achievable-perf2}).
\end{theorem}

\section{EXTENSION TO MODELS WITH NON-UNIFORM CONTROL COST}
\label{s:ExtensionModelsWithControlCost}

%[MATHEMATICAL FORMULATION WITH NON-UNIFORM CONTROL COST]
In the last section, we sought the test that minimizes the expected stopping time as the probability of error goes to zero.  In this formulation, we are implicitly penalizing the observation at each time equally regardless of the selected control signal at that time.  Consequently, a natural question that arises is how the adoption of an unequal control cost function affects the structure of the optimal test and its asymptotic performance.  To answer this question, we now consider an arbitrary cost function $\,c: \mathcal{U} \rightarrow \mathbb{R}^{+}.$  The goal will now be to minimize $\,\mathbb{E}_i \left [ \sum\limits_{k=1}^N c \left( U_k \right ) \right ]\,$ rather than $\,\mathbb{E}_i \left [ N \right ]\,$ as in 
the \tdblue{last} section.

\subsection{\tdblue{Proposed Test} and Its \tdblue{Asymptotic Optimality}}
\label{ss:TestControlCost}
Our test in the current setup is almost exactly the same as in the previous section except that (\ref{eqn-opt-exploit-control-nocost}) is changed to
\begin{equation}
q_i^* \left( \cdot \vert \cdot \right) 
\ =\ \mathop{\mbox{argmax}}_{q \left( \cdot \vert \cdot \right)}
~\frac{\min\limits_{j \neq i}
~\sum\limits_{\tilde{y} \in \mathcal{Y}, u \in \mathcal{U}}
\mu_i^q \left( \tilde{y} \right) q \left( u \vert \tilde{y} \right)
D \left( p_i^u \left( \cdot \| \tilde{y} \right) \| p_j^u \left( \cdot \| \tilde{y} \right) \right)}
{\sum\limits_{\tilde{y} \in \mathcal{Y}, u \in \mathcal{U}} \mu_i^q \left( \tilde{y} \right) q \left( u \vert \tilde{y} \right) c(u)}.
\label{eqn-opt-exploit-control-cost}
\end{equation}
It is interesting to note that the control cost function only affects the causal control policy and not the stopping rule (\ref{eqn-stopping-nocost}) and the final decision rule.  A heuristic explanation of this fact will be given in Remark \ref{remark-cost}.

\begin{theorem}[\tdblue{Asymptotic Optimality}]  
\label{thm-optasymperf-cost}
\tdblue{Any sequence of tests with vanishing maximal probability of error, $P_{ {\scriptsize \mathrm{max}}} \ \rightarrow\ 0,$ satisfies for every $i = 0, \ldots, M-1\,$ that 
\begin{eqnarray}
\mathbb{E}_i \left [ \sum\limits_{k=1}^N c \left( U_k \right ) \right ] &\geq& 
\frac{-\log{\left( 
R_i
\right)}}
{
\max\limits_{q \left( \cdot \vert \cdot \right)}
~\frac{\min\limits_{j \neq i}
~\sum\limits_{\tilde{y} \in \mathcal{Y}, u \in \mathcal{U}}
\mu_i^q \left( \tilde{y} \right) q \left( u \vert \tilde{y} \right)
D \left( p_i^u \left( \cdot \| \tilde{y} \right) \| p_j^u \left( \cdot \| \tilde{y} \right) \right)}
{\sum\limits_{\tilde{y} \in \mathcal{Y}, u \in \mathcal{U}} \mu_i^q \left( \tilde{y} \right) q \left( u \vert \tilde{y} \right) c(u)}
} 
\left( 1+ o(1) \right)
\label{eqn-converse-cost1}		\\
&\geq&
\frac{-\log{\left( 
P_{ {\scriptsize \mathrm{max}}}
\right)}}
{
\max\limits_{q \left( \cdot \vert \cdot \right)}
~\frac{\min\limits_{j \neq i}
~\sum\limits_{\tilde{y} \in \mathcal{Y}, u \in \mathcal{U}}
\mu_i^q \left( \tilde{y} \right) q \left( u \vert \tilde{y} \right)
D \left( p_i^u \left( \cdot \| \tilde{y} \right) \| p_j^u \left( \cdot \| \tilde{y} \right) \right)}
{\sum\limits_{\tilde{y} \in \mathcal{Y}, u \in \mathcal{U}} \mu_i^q \left( \tilde{y} \right) q \left( u \vert \tilde{y} \right) c(u)}
} 
\left( 1+ o(1) \right).
\label{eqn-converse-cost2}
\end{eqnarray}}

\tdblue{For the test in (\ref{eqn-opt-control-exploit}), (\ref{eqn-opt-exploit-control-cost}), (\ref{eqn-sparse-exploration}), (\ref{eqn-stopping-nocost}), it holds that
\begin{equation}
P_{ {\scriptsize \mathrm{max}}}^*
\ \leq\ \frac{M-1}{T}.
\end{equation}
Furthermore, the test is asymptotically optimal.  In particular, as the threshold $\,T\,$ in (\ref{eqn-stopping-nocost}) approaches infinity, the test in (\ref{eqn-opt-control-exploit}), (\ref{eqn-opt-exploit-control-cost}), (\ref{eqn-sparse-exploration}), (\ref{eqn-stopping-nocost}) yields for each $\,i = 0, \ldots, M-1\,$ that
\begin{eqnarray}
\mathbb{E}_i \left [ \sum\limits_{k=1}^{N^*} c \left( U_k \right ) \right ]
&=&
\frac{-\log{\left( 
P_{ {\scriptsize \mathrm{max}}}^*
\right)}}
{
\max\limits_{q \left( \cdot \vert \cdot \right)}
~\frac{\min\limits_{j \neq i}
~\sum\limits_{\tilde{y} \in \mathcal{Y}, u \in \mathcal{U}}
\mu_i^q \left( \tilde{y} \right) q \left( u \vert \tilde{y} \right)
D \left( p_i^u \left( \cdot \| \tilde{y} \right) \| p_j^u \left( \cdot \| \tilde{y} \right) \right)}
{\sum\limits_{\tilde{y} \in \mathcal{Y}, u \in \mathcal{U}} \mu_i^q \left( \tilde{y} \right) q \left( u \vert \tilde{y} \right) c(u)}
}
\left( 1+ o(1) \right)
\label{eqn-achievable-perf-cost1}	\\
&=&
\left( \inf_{\left( N, \delta \right) \,\in\, \mathbb{C} \left( P_{ {\scriptsize \mathrm{max}}}^* \right)} 
\mathbb{E}_i \left [ \sum\limits_{k=1}^{N} c \left( U_k \right ) \right ]
\right)
(1+o(1)) \nonumber \\
&=& 
\frac{-\log{\left( 
R_i^*
\right)}}
{
\max\limits_{q \left( \cdot \vert \cdot \right)}
~\frac{\min\limits_{j \neq i}
~\sum\limits_{\tilde{y} \in \mathcal{Y}, u \in \mathcal{U}}
\mu_i^q \left( \tilde{y} \right) q \left( u \vert \tilde{y} \right)
D \left( p_i^u \left( \cdot \| \tilde{y} \right) \| p_j^u \left( \cdot \| \tilde{y} \right) \right)}
{\sum\limits_{\tilde{y} \in \mathcal{Y}, u \in \mathcal{U}} \mu_i^q \left( \tilde{y} \right) q \left( u \vert \tilde{y} \right) c(u)}
}
\left( 1+ o(1) \right)
\label{eqn-achievable-perf-cost2} \\
&=& \left( \inf_{\left( N, \delta \right) \,\in\,  \mathbb{C}' \left( R_0^*, \ldots, R_{M-1}^* \right)} 
\mathbb{E}_i \left [ \sum\limits_{k=1}^{N} c \left( U_k \right ) \right ]  \right)
(1+o(1)).
\nonumber
\end{eqnarray}}
\end{theorem}
\begin{remark}
\label{remark-cost}
A heuristic explanation of \tdblue{the new test} accounting for the varying control cost function can also be given \tdblue{similarly to} Remark \ref{remark-nocost}.  Since we now seek to minimize $\,\mathbb{E}_i \left [ \sum\limits_{k=1}^N c \left( U_k \right ) \right ]\,$ instead of the $\,\mathbb{E}_i \left [ N \right ],$ we can reinterpret (\ref{eqn-stopping-nocost}) as
\begin{equation}
N \ \triangleq\ \min\limits_{n \geq 1}~
	~e^{ 
	\left ( \sum\limits_{k=1}^n c \left( U_k \right ) \right )
	\left(
	\frac{
	~\min\limits_{j \neq \hat{i}_n}~\frac{1}{n}
	\sum\limits_{k=1}^n 
	\log{\left( 
		\frac{p_{\hat{i}_n}^{U_k} \left(Y_k \vert Y_{k-1} \right)}
		{p_{j}^{U_k} \left(Y_k \vert Y_{k-1} \right)}
	\right) }
	}
	{
		\frac{1}{n} \sum\limits_{k=1}^n c \left( U_k \right )
	}
	\right )
	}
	~>~T.
	\label{eqn-stopping-cost-alt}
\end{equation}
The new causal control policy in (\ref{eqn-opt-control-exploit}), (\ref{eqn-opt-exploit-control-cost}), (\ref{eqn-sparse-exploration}) will now be the optimal self-tuning control policy maximizing the limiting ``cost-inferential'' reward
\begin{equation}
	\lim\limits_{n \rightarrow \infty}~
	\frac{
	~\min\limits_{j \neq i}~\frac{1}{n}
	\sum\limits_{k=1}^n 
	\log{\left( 
		\frac{p_{i}^{U_k} \left(Y_k \vert Y_{k-1} \right)}
		{p_{j}^{U_k} \left(Y_k \vert Y_{k-1} \right)}
	\right) }
	}
	{
		\frac{1}{n} \sum\limits_{k=1}^n c \left( U_k \right )
	},
\label{eqn-costed-inferential-rewards}
\end{equation}
and the optimal limiting value of this reward is what appears as the coefficient in the characterization of the optimal asymptotic performance in Theorem \ref{thm-optasymperf-cost}.
\end{remark}

\begin{remark}
We can also modify the test in Section \ref{ss:TestControlCost}
 to meet distinctly predefined constraints on the risks by resorting to the stopping rule (\ref{eqn-stopping-hardriskconstraints}).  As all the risks vanish while satisfying the assumption in Theorem \ref{thm-optasymperf-nocost-hardrisk}, the resulting test is also \tdblue{asymptotically optimal,} i.e., it satisfies (\ref{eqn-achievable-perf-cost1}) and (\ref{eqn-achievable-perf-cost2}).
\end{remark}

\begin{example} 
 \label{ex-1}
%Example with controlled memoryless observations and non-uniform control cost that shows that the cost function can substantially affects the optimal control policy.
This example shows that the varying control cost can substantially affect the optimal causal control policy.  This example (without the cost function) is taken from \cite*{niti-atia-veer-ieeetac-2013} for a simple model for controlled observations in which each $\,p_i^u \left( \cdot \vert \tilde{y} \right),\ i = 0, \ldots, M-1,$ is independent of $\,\tilde{y}.$ There are three hypotheses, and $\,\mathcal{Y} = \left \{0, 1 \right \},\ \mathcal{U} = 
\tdblue{\left \{\alpha, \beta, \gamma \right \}}.$  \tdblue{Consider a pair} of reciprocal distributions $\,p, \bar{p}\,$ on $\,\mathcal{Y},$ where $\,p(0) = \bar{p}(1) = \epsilon,\ 0 < \epsilon < 1$.  We let $\,p_i^u = p,$ for $\,\tdblue{i=0, u=\alpha;\ i=1, u=\beta;\ i=2, u=\gamma},$ and $\,p_i^u = \bar{p}\,$ for all other $\,(i,u)$.

It is not hard to see that when all control signals are equally costly, it holds that
\begin{eqnarray}
	q_0^*
&=& \mathop{\mathrm{argmax}}_{q \left( \cdot \right)}
~\min\limits_{j = 1, 2}
~\sum\limits_{u = \tdblue{\alpha, \beta, \gamma}}
q \left( u \right)
D \left( p_0^u  \| p_j^u  \right)
\nonumber \\
&=& 
\mathop{\mathrm{argmax}}_{q \left( \cdot \right)}~
\min \left \{
	  q \left( \tdblue{\alpha} \right) D \left( p \| \bar{p} \right) 
	  + q \left( \tdblue{\beta} \right) D \left( \bar{p} \| p \right), 
	  q \left( \tdblue{\alpha} \right) D \left( p \| \bar{p} \right)
	  + q \left( \tdblue{\gamma} \right) D \left( \bar{p} \| p \right)
\right \}	\nonumber \\
&=& 
\mathop{\mathrm{argmax}}_{q \left( \cdot \right)}~
\left(q \left( \tdblue{\alpha} \right) 
+ \min \left( q \left( \tdblue{\beta} \right), q \left( \tdblue{\gamma} \right) \right) \right)
\left( 1- 2 \epsilon \right) \log{ \left( \frac{1-\epsilon}{\epsilon} \right) }	
\nonumber \\
&=&  \mathbb{I}\left \{ u = \tdblue{\alpha} \right \}.	\nonumber 
\end{eqnarray}
By symmetry, we also get that $\,q_1^* = 
\mathbb{I}\left \{ u = \tdblue{\beta} \right \},$ 
and $\,q_2^* = \mathbb{I}\left \{ u = \tdblue{\gamma} \right \}.$

On the other hand, when the control cost function $\,c\,$ is not uniform, it can easily be shown that if $\,\tdblue{c(\alpha) < c(\beta)+c(\gamma)},$ then 
\begin{eqnarray}
q_0^{**}
&=& \mathop{\mathrm{argmax}}_{q \left( \cdot \right)}
~\frac{\min\limits_{j = 1, 2}
~\sum\limits_{u = \tdblue{\alpha, \beta, \gamma}}
q \left( u \right)
D \left( p_0^u  \| p_j^u  \right)}
{
\sum\limits_{u = \tdblue{\alpha, \beta, \gamma}}
q \left( u \right) c(u)
}	\nonumber \\
&=&  \mathbb{I}\left \{ u = \alpha \right \}.	\nonumber 
\end{eqnarray}
In other words, if the control signal $\,\tdblue{\alpha}\,$ is not too costly 
$\tdblue{(c(\alpha) < c(\beta)+c(\gamma))}$, it is still the optimal choice when the ML estimate of the hypothesis is $\,0\,$ (cf. (\ref{eqn-opt-exploit-control-cost})); otherwise, it will be better to start sampling from the less expensive control signals $\,\tdblue{\beta}\,$ and $\,\tdblue{\gamma}$.
\end{example}

\section{DISCUSSION}

Note that \tdblue{the optimal causal control policy} in (\ref{eqn-opt-control-exploit}), (\ref{eqn-opt-exploit-control-cost}), (\ref{eqn-sparse-exploration}) entails switching infinitely often between the ``exploiting'' control policy in (\ref{eqn-opt-control-exploit}), (\ref{eqn-opt-exploit-control-cost}) and the ``exploring'' control policy in (\ref{eqn-sparse-exploration}).  Such switching recurrences may not be practical for certain applications.  A control policy \tdgreen{that} completely separates the exploitation phase from the initial exploration phase with just one switch was proposed in \cite*{kief-sack-amstat-1963}[Section 5] for the simple observation model for controlled sensing for multihypothesis testing (as studied in \cite*{cher-amstat-1959, bess-tech-repo-1960, niti-atia-veer-ieeetac-2013, nagh-javi-astat-2013}).  The idea behind such a control policy should be directly extendable to the \tdblue{current controlled Markovian observation model.}  Another idea that can potentially be used to reduce or eliminate such switching is to select the instantaneous control signal based on a suitably biased version of the maximum likelihood estimate of the hypothesis, as in the self-tuning control policy proposed in \cite*{kuma-beck-ieeetac-1982}.  However, the biasing in \cite*{kuma-beck-ieeetac-1982} is based on the value function that is the solution of the dynamic programming equation associated with a {\em single} reward function.  On the other hand, \tdblue{the inferential reward} is a (non-linear) function of time averages with respect to {\em multiple} functions, namely $\,D \left( p_i^u \left( \cdot \vert \tilde{y} \right) \| p_j^u \left( \cdot \vert \tilde{y} \right) \right),\ j \neq i,$ and $\,c (u);$ hence, there is no suitable dynamic programming equation associated with \tdblue{the inferential reward maximization.}  Consequently, extending the idea of the biased maximum likelihood method in \cite*{kuma-beck-ieeetac-1982} to obtain an optimal control policy for \tdblue{the controlled sensing problem} remains a challenging future research problem.

\tdblue{The design of the asymptotically optimal test relied critically on the complete knowledge of the observation model.  In the future, it will be interesting to study how to design a good test when the knowledge of the distributional model for the controlled observations is incomplete. We also hope to explore whether the two-pronged approach of combining tools and principles from decision theory and stochastic control can be employed to study controlled sensing for solving other inference problems such as parameter estimation, and learning-based classification.}

%My talk at IWSM was received very well -- this was the right audience for the talk. There was quite a bit of discussion after the talk during which Yajun Mei brought to my attention a paper by Keifer and Sacks (see attached) in which the design of experiments problem is discussed in Section 5. This paper refers to Bessler, by the way. An interesting alternative to Chernoff's test is mentioned in this paper, in which the control is chosen based on a first block of observations, and then essentially fixed (with occasional exploration) as decision-making about the hypothesis is done. One reason for fixing the control is that in some applications (clinical trials) there is a cost associated with switching the experiment, which Chernoff did not consider. 

\section*{ACKNOWLEDGEMENTS}
%Thanks Ya Jun Mei for pointer to Kiefer and Sacks, and possibly Prof. Kumar if the work of Becker and Kumar would be relevant.
\tmred{We thank the associate editor and the referees for carefully reading the paper and providing useful and critical feedback.}  We also thank Professors P. R. Kumar and Ya Jun Mei for their pointers to \cite*{kuma-beck-ieeetac-1982} and \cite*{kief-sack-amstat-1963}, respectively.

This research was partially supported by the Air Force Office of Scientific Research (AFOSR) under Grant FA9550-10-1-0458 through the University of Illinois at Urbana-Champaign, by the U.S. Defense Threat Reduction Agency through subcontract 147755 at the University of Illinois from prime reward HDTRA1-10-1-0086, and by the National Science Foundation under Grant NSF CCF 11-11342.

\section*{APPENDIX TO SECTIONS \ref{s:AsymptoticallyOptimalTest}, \ref{s:ExtensionModelsWithControlCost}}

It is clear that Theorem \ref{thm-optasymperf-nocost} is a particularization of Theorem \ref{thm-optasymperf-cost} when the cost function $\,c(u) = 1,$ for every $\,u \in \mathcal{U}$.  Consequently, it suffices just to prove Theorem 
\ref{thm-optasymperf-cost}, and that the test in Section \ref{ss:TestControlCost} but with the stopping rule (\ref{eqn-stopping-hardriskconstraints}) satisfies the distinctly predefined constraints on the risks (\ref{eqn-hard-risk-constraints}) and achieves the \tdblue{asymptotic optimality} in (\ref{eqn-achievable-perf-cost1}) and (\ref{eqn-achievable-perf-cost2}).

For $\,i = 0, \ldots, M-1,$ let 
\begin{equation}
	d^*_i \ \triangleq\ 
	\max\limits_{q \left( \cdot \vert \cdot \right)}
~\frac{\min\limits_{j \neq i}
~\sum\limits_{\tilde{y} \in \mathcal{Y},\ u \in \mathcal{U}}
\mu_i^q \left( \tilde{y} \right) q \left( u \vert \tilde{y} \right)
D \left( 
\tmred{p_i^u \left( \cdot \vert \tilde{y} \right) \| p_j^u \left( \cdot \vert \tilde{y} \right)} 
\right)}
{\sum\limits_{\tilde{y} \in \mathcal{Y},\ u \in \mathcal{U}} 
\mu_i^q \left( \tilde{y} \right) q \left( u \vert \tilde{y} \right) c(u)}.
\label{eqn-target-opt-coeffs}	
\end{equation}

\vspace{0.2in}
\noindent
{\em Proof of the Converse Assertion in Theorem \ref{thm-optasymperf-cost}}
\vspace{0.1in}

We start with the proof of the \tdblue{converse assertion (\ref{eqn-converse-cost1})} which relies on the following two lemmas.
\begin{lemma}  
\label{proof-Thm5.1-lem1}
For any sequence of tests satisfying $\,P_{ {\scriptsize \mathrm{max}}} \rightarrow 0,$ and any small $\,\epsilon > 0,$ it holds for every $\,i \neq j,\ i, j \in \left \{ 0, \ldots, M-1 \right \},$ that
\begin{equation}
	\mathbb{P}_i \left \{
		\log{\left(\frac{p_i \left( Y^N, U^N \right)}{p_j \left( Y^N, U^N \right)} \right)}
		\leq - \left( 1 - \epsilon \right) \log{R_i}
	\right \}	\ \rightarrow\ 0.
\label{eqn-proof-Thm5.1-lem1}
\end{equation}
\end{lemma}

\begin{lemma}  
\label{proof-Thm5.1-lem2}
For any causal control policy, it holds for every $\,i = 0, \ldots, M-1\,$ and any small $\,\epsilon > 0\,$ that
\begin{equation}
%\lim\limits_{n \rightarrow \infty}
\mathbb{P}_i \left \{
~\min\limits_{j \neq i}~
\log{\left( \frac{p_i \left( Y^n, U^n \right)}{p_j \left( Y^n, U^n \right)} \right)}
\geq \left( \sum\limits_{k=1}^n c \left( U_k \right) \right) \left( 
	d^*_i + \epsilon 
\right)
\right \}  \ \leq\ O \left( n^{-2} \right).
\label{eqn-proof-Thm5.1-lem2}
\end{equation}
\end{lemma}

\tdblue{Lemma \ref{proof-Thm5.1-lem1} follows from observing that}
\tdblue{\begin{align}
	\mathbb{P}_i \left \{
		\log{\left(\frac{p_i \left( Y^N, U^N \right)}{p_j \left( Y^N, U^N \right)} \right)}
		\leq - \left( 1 - \epsilon \right) \log{R_i}
	\right \}		\hspace{2.8in}		\nonumber \\
	\ \leq\ 
	\mathbb{P}_i \left \{
		\log{\left(\frac{p_i \left( Y^N, U^N \right)}{p_j \left( Y^N, U^N \right)} \right)}
		\leq - \left( 1 - \epsilon \right) 
		\log{\left( 
			\sum\limits_{j \neq i} \frac{1}{M-1} \mathbb{P}_j \left \{ \delta = i \right \} 
		\right)}
	\right \} \ \rightarrow\ 0.
\label{eqn-proof-Thm5.1-lem1-pf}
\end{align}}
\tdblue{The first identity in (\ref{eqn-proof-Thm5.1-lem1-pf}) follows directly from the definition of $\,R_i,\ i = 0, \ldots, M-1,$ where as the second identity follows as in \cite*{niti-atia-veer-ieeetac-2013}[Lemma 3]} by working instead on the joint distribution of observations and control signals $\,p_i \left( y^n, u^n \right),\ i = 0, \ldots, M-1$ (cf. \cite*{niti-atia-veer-ieeetac-2013}[Lemma 3]).  The main effort in the converse proof will be in proving Lemma \ref{proof-Thm5.1-lem2}; however, we first show how \tdblue{the converse assertion} follows from Lemmas \ref{proof-Thm5.1-lem1} and \ref{proof-Thm5.1-lem2}.
Let $\,c_{R_i} = \frac{- \left( 1 - \epsilon \right) \log{R_i}}{d^*_i + \epsilon},\ \underline{c} = \min\limits_{u \in \mathcal{U}} c(u),$ and $\underline{p} = 
\min\limits_{u \in \mathcal{U},\ y, \tilde{y} \in \mathcal{Y},\ i = 0, \ldots, M-1} p_i^u \left( y \vert \tilde{y} \right)$.  
Then, we get for every $\,i = 0, \ldots, M-1\,$ that
\begin{eqnarray}
%CHECK THIS STEP
\mathbb{P}_i \left \{ \sum\limits_{k=1}^N c \left( U_k \right) \leq c_{R_i} \right \}
&\leq&
\mathbb{P}_i \left \{ \sum\limits_{k=1}^N c \left( U_k \right) \leq c_{R_i}, 
\min\limits_{j \neq i}
\log{\left( \frac{p_i \left( Y^N, U^N \right)}{p_j \left( Y^N, U^N \right)} \right)} 
\geq - \left( 1- \epsilon \right) \log{R_i}
\right \}	\nonumber \\
& & +\  
\mathbb{P}_i \left \{ 
\min\limits_{j \neq i}
\log{\left( \frac{p_i \left( Y^N, U^N \right)}{p_j \left( Y^N, U^N \right)} \right)} 
\leq - \left( 1- \epsilon \right) \log{R_i}
\right \}	\nonumber \\
&\leq& \mathbb{P}_i \left \{  
\begin{array}{cc}
\min\limits_{j \neq i}
\log{\left( \frac{p_i \left( Y^N, U^N \right)}{p_j \left( Y^N, U^N \right)} \right)} 
\geq \left( \sum\limits_{k=1}^N c \left ( U_k \right) \right) \left( d^*_i + \epsilon \right); \\
\frac{- \left( 1 - \epsilon \right) \log{R_i}}{\log{\left( \frac{1}{\underline{p}} \right)}}
\leq
N \leq \frac{- \left( 1 - \epsilon \right) \log{R_i}}{\left( d^*_i + \epsilon \right) \underline{c}}
\end{array}
\right \}	\nonumber \\
& & +\  
\mathbb{P}_i \left \{ 
\min\limits_{j \neq i}
\log{\left( \frac{p_i \left( Y^N, U^N \right)}{p_j \left( Y^N, U^N \right)} \right)} 
\leq - \left( 1- \epsilon \right) \log{R_i}
\right \}.	
\label{eqn-proof-Thm5.1-conv}
\end{eqnarray}
%\texttt{CHECK THIS STEP INVOKING LEMMA \ref{proof-Thm5.1-lem2}}\\
Lemma \ref{proof-Thm5.1-lem1} yields that the second term on the right of (\ref{eqn-proof-Thm5.1-conv}) vanishes to zero, while Lemma \ref{proof-Thm5.1-lem2} gives us that for every $\,\eta > 0,$ the limit of the first term (as $\,R_i \rightarrow 0$) on the right of (\ref{eqn-proof-Thm5.1-conv}) is  bounded above by a term of order $\,O \left( \frac{1}{ \sqrt{\eta}} \eta \right) = O \left( \sqrt{\eta} \right),$ yielding also that the second term is going to zero.  The converse assertion is now proved.
%Lemma \ref{proof-Thm5.1-lem2} and \ref{proof-Thm5.1-lem1} yield that both the first and the second terms on the right-side of (\ref{eqn-proof-Thm5.1-conv}) vanish to zero.  The converse proof now follows, as $\,\epsilon\,$ can be arbitrarily small.

We now prove Lemma \ref{proof-Thm5.1-lem2}.  To this end, we note that for any $\,0 \leq i \neq j \leq M-1,$
\begin{eqnarray}
	\log{\left( \frac{p_i \left( Y^n, U^n \right)}{p_j \left( Y^n, U^n \right) } \right)}
	&=&
	\sum\limits_{k = 1}^n
	~\log{\left( \frac{p_i^{U_k} \left( Y_k \vert Y_{k-1} \right)}
	{p_j^{U_k} \left( Y_k \vert Y_{k-1} \right) } \right)}.	\nonumber \\
	&=&
	\sum\limits_{k = 1}^n
	\left \{
	\log{\left( \frac{p_i^{U_k} \left( Y_k \vert Y_{k-1} \right)}
	{p_j^{U_k} \left( Y_k \vert Y_{k-1} \right) } \right)}
	~-~ D \left( p_i^{U_k} \left( \cdot \vert Y_{k-1} \right) \|
				  p_j^{U_k} \left( \cdot \vert Y_{k-1} \right) \right)
	\right \}	\nonumber \\
	& & +\ \sum\limits_{k = 1}^n
		D \left( p_i^{U_k} \left( \cdot \vert Y_{k-1} \right) \|
				  p_j^{U_k} \left( \cdot \vert Y_{k-1} \right) \right).	
	\label{eqn-PfLemm6.4-a}
\end{eqnarray}
For $\,j \neq i,$ with
%We now apply Claim 1 with $\,\mathcal{M} = \left \{ j,\ j \neq i \right \},$
\begin{eqnarray}
A_k^j &=& \left \{
	\log{\left( \frac{p_i^{U_k} \left( Y_k \vert Y_{k-1} \right)}
	{p_j^{U_k} \left( Y_k \vert Y_{k-1} \right) } \right)}
	- D \left( p_i^{U_k} \left( \cdot \vert Y_{k-1} \right) \|
				  p_j^{U_k} \left( \cdot \vert Y_{k-1} \right) \right)
	\right \},\nonumber \\ 
B_k^j &=& D \left( p_i^{U_k} \left( \cdot \vert Y_{k-1} \right) \|
				  p_j^{U_k} \left( \cdot \vert Y_{k-1} \right) \right),	\nonumber
\end{eqnarray}
$k = 1, \ldots, n,$ we get that
\begin{equation}
\mathbb{P}_i \left \{
%\max\limits_{1 \leq m \leq n}
~\min\limits_{j \neq i}~
\log{\left( \frac{p_i \left( Y^n, U^n \right)}{p_j \left( Y^n, U^n \right)} \right)}
\geq \left( \sum\limits_{k=1}^n c \left( U_k \right) \right) \left( 
	d^*_i + \epsilon 
\right)
\right \}	\nonumber
\hspace{1.3in}
\end{equation}
\begin{eqnarray}
&\leq& \sum\limits_{j \neq i}~
\mathbb{P}_i \left \{
%\max\limits_{1 \leq m \leq n}
\left( \sum\limits_{k=1}^n A_k^j \right)
\ \geq\ \left( \sum\limits_{k=1}^n c \left( U_k \right) \right) \frac{\epsilon}{2}
\right \}	\nonumber \\
& & +\ 
\mathbb{P}_i \left \{
\min\limits_{j \neq i} \left( \sum\limits_{k=1}^n B_k^j \right)
\ \geq\ \left( \sum\limits_{k=1}^n c \left( U_k \right) \right) 
\left( d_i^* +  \frac{\epsilon}{2} \right)
\right \}	\nonumber \\
&\leq& \sum\limits_{j \neq i}~
\mathbb{P}_i \left \{
%\max\limits_{1 \leq m \leq n}
\left( \sum\limits_{k=1}^n A_k^j \right)
\ \geq\ n \underline{c} \frac{\epsilon}{2}
\right \}	\nonumber \\
& & +\ 
\mathbb{P}_i \left \{
\min\limits_{j \neq i} \left( \sum\limits_{k=1}^n B_k^j \right)
\ \geq\ \left( \sum\limits_{k=1}^n c \left( U_k \right) \right) 
\left( d_i^* +  \frac{\epsilon}{2} \right)
\right \}.	
\label{eqn-PfLem6.4-a}
\end{eqnarray}
To finish the proof of Lemma \ref{proof-Thm5.1-lem2}, it suffices to prove that each probability on the right-side of (\ref{eqn-PfLem6.4-a}) vanishes exponentially in $\,n$.

We first prove that each of the summands in the first term on the right-side of (\ref{eqn-PfLem6.4-a}) vanishes exponentially rapidly.  To this end, we shall invoke the Chernoff bounding argument.  Note that for any $\,\epsilon > 0,$ and any $\,j \neq i,\ \tilde{y} \in \mathcal{Y},\ u \in \mathcal{U},$ 
\begin{equation}
\mathbb{E}_i \left [
- \log{\left( \frac{p_i^{u} \left( Y_k \vert \tilde{y} \right)}
	{p_j^{u} \left( Y_k \vert \tilde{y} \right) } \right)}
	+ D \left( p_i^{u} \left( \cdot \vert \tilde{y} \right) \|
				  p_j^{u} \left( \cdot \vert \tilde{y} \right) \right)
	+ \underline{c} \frac{\epsilon}{2}
\bigg \vert Y_{k-1} = \tilde{y}, U_k = u
\right ] \ =\ \underline{c} \frac{\epsilon}{2} \ >\ 0, 
\label{eqn-proof-Lemma6.5-b}
\end{equation}
and that the random variable inside the expectation in (\ref{eqn-proof-Lemma6.5-b}) yields a finite value of the conditional moment generating function for the parameter range $\,-1 \leq t \leq 0\,$ (cf. (\ref{eqn-assump-pos-transition})).  Since, $\,\mathcal{U}\,$ and $\,\mathcal{Y}\,$ are both finite, there are $\,t^* (\epsilon) < 0\,$ and $\,b(\epsilon) > 0\,$ such that all the conditional moment generating functions (among all possible $\,\tilde{y}\,$ and $\,u$) evaluated at $\,t^* (\epsilon)\,$ are uniformly bounded by $\,e^{-b} < 1$.  Hence,
\begin{equation}
\mathbb{P}_i \left \{
\sum\limits_{k=1}^n
\left \{
	- \log{\left( \frac{p_i^{U_k} \left( Y_k \vert Y_{k-1} \right)}
	{p_j^{U_k} \left( Y_k \vert Y_{k-1} \right) } \right)}
	+ D \left( p_i^{U_k} \left( \cdot \vert Y_{k-1} \right) \|
				  p_j^{U_k} \left( \cdot \vert Y_{k-1} \right) \right)
	+ \underline{c} \frac{\epsilon}{2} 
\right \}
\ \leq\ 0
\right \}
\hspace{1.0in}
\nonumber 
\end{equation}
\begin{eqnarray}
&\leq&
\mathbb{E}_i \left [ 
e^{ 
t\left( \epsilon \right) 
\left (
\sum\limits_{k=1}^n
\left \{
	- \log{\left( \frac{p_i^{U_k} \left( Y_k \vert Y_{k-1} \right)}
	{p_j^{U_k} \left( Y_k \vert Y_{k-1} \right) } \right)}
	+ D \left( p_i^{U_k} \left( \cdot \vert Y_{k-1} \right) \|
				  p_j^{U_k} \left( \cdot \vert Y_{k-1} \right) \right)
	+ \underline{c} \frac{\epsilon}{2} 
\right \}
\right )
}
\right ].\nonumber \\
&=& 
\mathbb{E}_i \left [ 
\mathbb{E}_i \left [ 
e^{ 
t\left( \epsilon \right) 
\left (
\sum\limits_{k=1}^n
\left \{
	- \log{\left( \frac{p_i^{U_k} \left( Y_k \vert Y_{k-1} \right)}
	{p_j^{U_k} \left( Y_k \vert Y_{k-1} \right) } \right)}
	+ D \left( p_i^{U_k} \left( \cdot \vert Y_{k-1} \right) \|
				  p_j^{U_k} \left( \cdot \vert Y_{k-1} \right) \right)
	+ \underline{c} \frac{\epsilon}{2} 
\right \}
\right )
}
\bigg \vert Y^{n-1}, U^n
\right ]
\right ]		\nonumber \\
&=& 
e^{-b}~
\mathbb{E}_i \left [ 
e^{ 
t\left( \epsilon \right) 
\left (
\sum\limits_{k=1}^{n-1}
\left \{
	- \log{\left( \frac{p_i^{U_k} \left( Y_k \vert Y_{k-1} \right)}
	{p_j^{U_k} \left( Y_k \vert Y_{k-1} \right) } \right)}
	+ D \left( p_i^{U_k} \left( \cdot \vert Y_{k-1} \right) \|
				  p_j^{U_k} \left( \cdot \vert Y_{k-1} \right) \right)
	+ \underline{c} \frac{\epsilon}{2} 
\right \}
\right )
}
\right ]	.	\label{eqn-proof-Lemma6.5-b-2}
\end{eqnarray}
Successive uses of the smoothing property of conditional expectation then yield that each of the summands in the first term on the right-side of (\ref{eqn-PfLem6.4-a}) converges to zero exponentially.  This proof step is  similar to the one in the proof of \cite*{cher-amstat-1959}[Equation (5.10)].

%The proof that the first term on the right-side of (\ref{eqn-PfLem6.4-a}) converges to zero follows from a standard martingale convergence argument, since for each $\,j \neq i,$ the sequence 
%\begin{equation}
%S_m = \sum\limits_{k=1}^m 
%	\left \{
%	\log{\left( \frac{p_i^{U_k} \left( Y_k \vert Y_{k-1} \right)}
%	{p_j^{U_k} \left( Y_k \vert Y_{k-1} \right) } \right)}
%	- D \left( p_i^{U_k} \left( \cdot \vert Y_{k-1} \right) \|
%				  p_j^{U_k} \left( \cdot \vert Y_{k-1} \right) \right)
%	 \right \}	\nonumber
%\end{equation}
%is a martingale (admissible for the filtration generated by the observation and control sequences) with finite incremental conditional means and conditional variances (cf. (\ref{eqn-assump-pos-transition})).  
%This argument is very similar to that in the proof of \cite*{cher-amstat-1959}[Lemma 5].

It is now left only to prove that 
\begin{equation}
%\lim\limits_{n \rightarrow \infty}
\mathbb{P}_i \left \{
\min\limits_{j \neq i}~\sum\limits_{k=1}^n 
D \left( p_i^{U_k} \left( \cdot \vert Y_{k-1} \right) \|
				  p_j^{U_k} \left( \cdot \vert Y_{k-1} \right) \right)
\ \geq\ \left( \sum\limits_{k=1}^n c \left( U_k \right) \right) 
\left( d_i^* +  \frac{\epsilon}{2} \right)
\right \}	
%\ =\ 0.
\label{eqn-proof-Lemma7.4-b-2-aa}
\end{equation}
vanishes exponentially rapidly.
%Note that it suffices to prove that for any causal control policy,
%\begin{equation}
%\limsup\limits_{n \rightarrow \infty}
%\frac{
%\min\limits_{j \neq i}~
%\frac{1}{n} \sum\limits_{k=1}^n 
%D \left( p_i^{U_k} \left( \cdot \vert Y_{k-1} \right) \|
%				  p_j^{U_k} \left( \cdot \vert Y_{k-1} \right) \right)
%}
%{
%\frac{1}{n} \sum\limits_{k=1}^n c \left( U_k \right)
%}
%\ \ \leq\ \ d_i^*\ \ \mbox{a.s.}
%\label{eqn-PfLem6.4-c}
%\end{equation}
To this end, we resort to the convex analytic method used widely in Markov Decision Processes (see, e.g., \cite*{bork-cmc-book-1991}) which entails looking at the joint empirical distribution of observation and control sequences $\,T_n \left( \tilde{y}, u \right), \tilde{y} \in \mathcal{Y}, u \in \mathcal{U},$ defined as
\begin{equation}
	T_n \left( \tilde{y}, u \right) \ =\
	\frac{1}{n} \sum\limits_{k=1}^n  \mathbb{I} \left \{ Y_{k-1} = \tilde{y}, U_k = u \right \}.
	\label{eqn-def-randomjointemp} 
\end{equation}
%MODIFIED TILL HERE (1)
First, note that it follows from (\ref{eqn-target-opt-coeffs}) 
%and Assumption (\ref{eqn-assump-pos-transition}) 
upon identifying $\,t\left( \tilde{y}, u \right) = \mu_i^q \left( \tilde{y} \right)
q \left( u \vert \tilde{y} \right)\,$
that 
\begin{equation}
d_i^* \ =\
\max\limits_{
t \left( \tilde{y}, u \right):\
t \left( y \right) ~=~ 
\sum\limits_{\tilde{y}, u} t \left( \tilde{y}, u \right)
p_i^u \left( y \vert \tilde{y} \right)
}~
\frac{
\sum\limits_{\tilde{y} \in \mathcal{Y}, u \in \mathcal{U}}
t\left( \tilde{y}, u \right) 
D \left( p_i^u \left( \cdot \vert \tilde{y} \right) \|
p_j^u \left( \cdot \vert \tilde{y} \right)
\right)
}
{
\sum\limits_{u \in \mathcal{U}}
t(u) c(u)
}.
\label{eqn-proof-Lemma7.4-b-2-a}
\end{equation}
%where $\,\mu_i^{q_i^*}\,$ is the unique distribution satisfying for every $\,y \in \mathcal{Y}\,$ that
%\begin{equation}
%\mu_i^{q_i^*} \left( y \right)  
%= \sum\limits_{\tilde{y}\in \mathcal{Y}}
%p_i^{q^*_i} \left( y \vert \tilde{y} \right)~\mu_i^{q_i^*} \left( \tilde{y} \right).
%\label{eqn-proof-Lemma6.5-f-2}
%\end{equation}
Next, since $\, c \left( u \right ) \geq \underline{c} > 0,$ we get that 
the quantity 
$\,\frac{
\sum\limits_{\tilde{y} \in \mathcal{Y}, u \in \mathcal{U}}
t\left( \tilde{y}, u \right) 
D \left( p_i^u \left( \cdot \vert \tilde{y} \right) \|
p_j^u \left( \cdot \vert \tilde{y} \right)
\right)
}
{
\sum\limits_{u \in \mathcal{U}}
t(u) c(u)
}\,
$ 
is a continuous function of the joint distribution $\,t$.  By this continuity, 
it holds that there exists a sufficiently small $\,\eta\,$ such that
\begin{equation}
d_i^* + \frac{\epsilon}{2}  \ \ >\ \ \ 
\max\limits_{
t:\ 
\bigg \vert 
t \left( \tilde{y}, u \right):\
t \left( y \right) -
\sum\limits_{\tilde{y}, u} t \left( \tilde{y}, u \right)
p_i^u \left( y \vert \tilde{y} \right)
\bigg \vert \ 
\leq\  \eta,\ \forall y \in \mathcal{Y}
} \ \ \  
\frac{
t\left( \tilde{y}, u \right) 
D \left( p_i^u \left( \cdot \vert \tilde{y} \right) \|
p_j^u \left( \cdot \vert \tilde{y} \right)
\right)
}
{
\sum\limits_{u \in \mathcal{U}}
t(u) c(u)
}.
\label{eqn-proof-Lemma7.4-b-2-b}
\end{equation}
%Now we look at the empirical distribution
%\begin{equation}
%	T \left( \tilde{y} \right) \ =\
%	\frac{1}{n - \lfloor \epsilon' n \rfloor + 1} 
%	\sum\limits_{k={\lfloor \epsilon' n \rfloor}}^n  \mathbb{I} \left \{ Y_{k-1} = \tilde{y} \right \}.
%	\nonumber	
%\end{equation}
Then, if follows from (\ref{eqn-proof-Lemma7.4-b-2-b}) that
\begin{equation}
\mathbb{P}_i \left \{
\min\limits_{j \neq i}~\sum\limits_{k=1}^n 
D \left( p_i^{U_k} \left( \cdot \vert Y_{k-1} \right) \|
				  p_j^{U_k} \left( \cdot \vert Y_{k-1} \right) \right)
\ \geq\ \left( \sum\limits_{k=1}^n c \left( U_k \right) \right) 
\left( d_i^* +  \frac{\epsilon}{2} \right)
\right \}	
	\hspace{2.0in}
\nonumber
\end{equation}
\begin{eqnarray} 
&=&
\mathbb{P}_i 
\left \{
\frac{\sum\limits_{\tilde{y} \in \mathcal{Y}, u \in \mathcal{U}} 
T_n \left( \tilde{y}, u \right)
D \left( p_i^{u} \left( \cdot \vert \tilde{y} \right) \|
				  p_j^{u} \left( \cdot \vert \tilde{y} \right) \right)
}
{\sum\limits_{u \in \mathcal{U}} T_n \left( u \right) c \left( u \right)}
\geq d_i^* + \frac{\epsilon}{2}
\right \} 	\nonumber \\
&\leq&
\mathbb{P}_i 
\left \{
\bigg \vert 
T_n \left( y \right) -
\sum\limits_{\tilde{y} \in \mathcal{Y}, u \in \mathcal{U}} 
T_n \left( \tilde{y}, u \right)
p_i^u \left( y \vert \tilde{y} \right)
\bigg \vert 
\ >\  \eta,\ \exists y \in \mathcal{Y}
\right \}.
\label{eqn-proof-Lemma7.4-b-2-c}
\end{eqnarray}
The same Chernoff bounding argument applied to 
\begin{equation}
\sum\limits_{k=1}^n  
\left(
\mathbb{I} \left \{ Y_{k} = y \right \}
- \mathbb{E}_i \left [
\mathbb{I} \left \{ Y_{k} = y \right \}
\big \vert Y^{k-1}, U^k
\right ]
\right)
\hspace{2.0in}
\nonumber 
\end{equation}
\begin{equation}
\ \ \ = 
\sum\limits_{k=1}^n  
\left(
\mathbb{I} \left \{ Y_{k} = y \right \}
- \sum\limits_{\tilde{y} \in \mathcal{Y}, u \in \mathcal{U}}
p_i^u \left( y \vert \tilde{y} \right) 
~\mathbb{I} \left \{ Y_{k-1} = \tilde{y}, U_k = u \right \}
\right)
\nonumber
\end{equation}
yields that the probability on the right-side of (\ref{eqn-proof-Lemma7.4-b-2-c}) vanishes exponentially in $\,n,$ thereby proving (\ref{eqn-proof-Lemma7.4-b-2-aa}) and Lemma \ref{proof-Thm5.1-lem2}.

\vspace{0.2in}
\noindent
{\em Proof of the Forward Assertion in Theorem \ref{thm-optasymperf-cost}}
\vspace{0.1in}

Next, we move on the the proof of the forward assertion (\ref{eqn-achievable-perf-cost1}) which relies on the following Lemmas.

\begin{lemma}  
\label{proof-Thm5.1-lem3}
When the causal control policy (\ref{eqn-opt-control-exploit}), (\ref{eqn-opt-exploit-control-cost}), (\ref{eqn-sparse-exploration}) is applied perpetually, it holds for any $\, \gamma > 2,$ every $\,0 \leq i \neq j \leq M-1\,$ and any small $\,\epsilon > 0\,$ that
\begin{equation}
\lim\limits_{n \rightarrow \infty}
\mathbb{P}_i \left \{
\log{\left( \frac{p_i \left( Y^n, U^n \right)}{p_j \left( Y^n, U^n \right)} \right)}
\leq \left( \sum\limits_{k=1}^n c \left( U_k \right) \right) \left( 
	d^*_i - \epsilon 
\right)
\right \}  \ =\ O \left( n^{- \gamma} \right).
\label{eqn-proof-Thm5.1-lem3}
\end{equation}
\end{lemma}

\begin{lemma}  
\label{proof-Thm5.1-forward-lem4}  For any $\,\gamma > 2,$ when the parameter $\,a\,$ in (\ref{eqn-sparse-exploration}) is selected to be sufficiently close to 1, it holds for every $\,i = 0, \ldots, M-1,$ and any $\,\epsilon' > 0\,$ that the first time $\,\tilde{T}\,$ after which the ML estimate $\,\hat{i}\,$ always equals the true hypothesis $\,i,$ satisfies
\begin{equation}
\lim\limits_{n \rightarrow \infty}
\mathbb{P}_i \left \{
\tilde{T} \ \geq\ \epsilon' n
\right \}  \ =\ O \left( n^{- \gamma} \right).
\label{eqn-proof-Thm5.1-forward-lem4}
\end{equation}
\end{lemma}

First, we show that Lemmas \ref{proof-Thm5.1-lem3}, \ref{proof-Thm5.1-forward-lem4} will lead to (\ref{eqn-achievable-perf-cost1}).  This finishing step of the proof follows the ideas in the proof of \cite*{cher-amstat-1959}[Lemma 2].  To this end, for a fixed $\,i = 0, \ldots, M-1,$ and $\,j \neq i,$ we let $\,N_j\,$ be the smallest time for which $\,\frac{p_i \left( Y^n, U^n \right)}{p_j \left( Y^n, U^n \right)} > T,$ for all $\,n \geq N_j$.  Hence, we get that 
$\,\tdblue{N^*} \leq \max \left( \max\limits_{j \neq i} N_j, \tilde{T} \right)$.  Now for any $\,\epsilon > 0\,$ and $\,A > \frac{\log{T}}{\left( d_i^* - \epsilon \right)}\,$ with $\,\overline{c} = \max\limits_{u \in \mathcal{U}} c(u),$ we get from Lemmas \ref{proof-Thm5.1-forward-lem4}, \ref{proof-Thm5.1-lem3} (with the $\,\gamma\,$ therein being the same), respectively, that
\begin{equation}
\mathbb{P}_i \left \{ \sum\limits_{k=1}^{\tilde{T}} c \left( U_k \right) > A \right \}
\ \leq\ \mathbb{P}_i \left \{ \tilde{T} > \frac{A}{\overline{c}} \right \}
\ =\ O \left( \left( \frac{A}{\overline{c}} \right)^{-\gamma} \right),
\label{eqn-proof-Thm5.1-forward-a-0}
\end{equation}
and that
\begin{eqnarray}
\mathbb{P}_i \left \{ \sum\limits_{k=1}^{N_j} c \left( U_k \right) >  A \right \}
&\leq& 
\mathbb{P}_i \left \{ 
\sum\limits_{k=1}^{N_j - 1} c \left( U_k \right) > 
A - \overline{c};\ \left( N_j - 1 \right) \geq \lfloor \frac{A}{\overline{c}} \rfloor
\right \}	\nonumber \\
%\texttt{CHECK} \nonumber \\ 
&\leq& 
\mathbb{P}_i \left \{ 
\sum\limits_{k=1}^{N_j - 1} c \left( U_k \right) > 
\frac{\log{T}}{\left( d_i^* - \epsilon \right)} - \overline{c};\ 
\left( N_j - 1 \right) \geq \lfloor \frac{A}{\overline{c}} \rfloor
\right \}	\nonumber \\
&\leq&  
\sum\limits_{n = \lfloor \frac{A}{\overline{c}} \rfloor}^{\infty}
\mathbb{P}_i \left \{
\sum\limits_{k=1}^{n} c \left( U_k \right) > 
\frac{1}{\left( d_i^* - \frac{\epsilon}{2} \right)}
\log{\left( \frac{p_i \left( Y^n, U^n \right)}{p_j \left( Y^n, U^n \right)} \right)}
\right \} 	\nonumber \\
&=&
\sum\limits_{n = \lfloor \frac{A}{\overline{c}} \rfloor}^{\infty}
O \left( n^{-\gamma} \right)
\ =\ O \left( \left( \frac{A}{\overline{c}} \right)^{-\gamma +1} \right).
\label{eqn-proof-Thm5.1-forward-a}
\end{eqnarray}
Following from (\ref{eqn-proof-Thm5.1-forward-a-0}), (\ref{eqn-proof-Thm5.1-forward-a}), we get for any $\,\epsilon > 0\,$ that
\begin{eqnarray}
\mathbb{E}_i \left [ \sum\limits_{k=n}^{\tilde{T}} c \left( U_k \right) \right ], 
\mathbb{E}_i \left [ \sum\limits_{k=n}^{N_j} c \left( U_k \right) \right ]
&\leq&  \frac{\log{T}}{\left(d_i^* - \epsilon\right)}
\left \{ 1 + 
\int\limits_{\frac{\log{T}}{\left(d_i^* - \epsilon\right)}}^{\infty}
O\left( \left( \frac{A}{\overline{c}} \right)^{-\gamma+1} \right) dA 
\right \}
\nonumber \\
&=& \frac{\log{T}}{\left(d_i^* - \epsilon\right)}
\left \{ 1 + 
O\left( \left( \frac{\log{T}}{\overline{c}\left(d_i^* - \epsilon\right)} \right)^{-\gamma+2} \right)
\right \}	\nonumber \\
&=& \frac{\log{T}}{\left(d_i^* - \epsilon\right)}
\left \{ 1 + 
o\left( \log{T} \right)
\right \}.	\nonumber 
\end{eqnarray}
Consequently, we get from $\,\tdblue{N^*} \leq \max \left( \max\limits_{j \neq i} N_j, \tilde{T} \right)\,$ that
\begin{equation}
\mathbb{E}_i \left [ \sum\limits_{k=n}^{\tdblue{N^*}} c \left( U_k \right) \right ]
\ \leq\ \frac{\log{T}}{\left(d_i^* - \epsilon\right)} (1 + o(1)).
\label{eqn-proof-Thm5.1-forward-b}
\end{equation}
Next, for each $\,j \neq i,$ let $\,A_{n, j} = \left \{ \left( y^n, u^n \right), 
\tdblue{N^*} = n, \delta = j \right \}.$  Then, we get from (\ref{eqn-stopping-nocost}) that for each $\,j \neq i,\ \mathbb{P}_i \left \{ A_{n, j} \right \} \leq \frac{1}{T} \mathbb{P}_j \left \{ A_{n, j} \right \}$.  Hence, we get that 
$\mathbb{P}_i \left \{ \delta \neq i \right \} = \sum\limits_{j \neq i} \sum\limits_{n=1}^{\infty} \mathbb{P}_i \left \{ A_{n, j} \right \} \leq \sum\limits_{j \neq i} \frac{1}{T} \sum\limits_{n=1}^{\infty} \mathbb{P}_j \left \{ A_{n, j} \right \} = (M-1) \frac{1}{T},$ and, hence, $ \tdblue{P_{ {\scriptsize \mathrm{max}}}^* \leq \frac{(M-1)}{T}}$.  Combining this with (\ref{eqn-proof-Thm5.1-forward-b}) yields (\ref{eqn-achievable-perf-cost1}).

We now prove Lemma \ref{proof-Thm5.1-lem3}.  To this end, we first note that
\begin{equation}
\mathbb{P}_i \left \{
\log{\left( \frac{p_i \left( Y^n, U^n \right)}{p_j \left( Y^n, U^n \right)} \right)}
\leq \left( \sum\limits_{k=1}^n c \left( U_k \right) \right) \left( 
	d^*_i - \epsilon 
\right)
\right \}
\hspace{2.2in}
\end{equation}
\begin{eqnarray}
&\leq&
\mathbb{P}_i \left \{
\sum\limits_{k=1}^n
\left \{
	\log{\left( \frac{p_i^{U_k} \left( Y_k \vert Y_{k-1} \right)}
	{p_j^{U_k} \left( Y_k \vert Y_{k-1} \right) } \right)}
	- D \left( p_i^{U_k} \left( \cdot \vert Y_{k-1} \right) \|
				  p_j^{U_k} \left( \cdot \vert Y_{k-1} \right) \right)
\right \}
\leq -n\underline{c} \frac{\epsilon}{2} 
\right \}	\nonumber \\
& & +\  
\mathbb{P}_i \left \{
\sum\limits_{k=1}^n
\left \{
D \left( p_i^{U_k} \left( \cdot \vert Y_{k-1} \right) \|
				  p_j^{U_k} \left( \cdot \vert Y_{k-1} \right) \right)
		- d_i^* c \left( U_k \right)
\right \}
\leq -n\underline{c} \frac{\epsilon}{2} 
\right \}.
\label{eqn-proof-Lemma6.5-a}	
\end{eqnarray}
The proof that the probability of the first term on the right-side of (\ref{eqn-proof-Lemma6.5-a}) goes to zero exponentially in $\,n\,$ can be carried out by invoking same Chernoff bounding argument as the one leading to (\ref{eqn-proof-Lemma6.5-b-2}).

So now it is left to prove that for any $\,\tilde{\epsilon} > 0$,
\begin{equation}
\mathbb{P}_i 
\left \{
\sum\limits_{k=1}^n
\left (
D \left( p_i^{U_k} \left( \cdot \vert Y_{k-1} \right) \|
				  p_j^{U_k} \left( \cdot \vert Y_{k-1} \right) \right)		 
- 
c \left( U_k \right)
d_i^*
+ \tilde{\epsilon}
\right )
\ \leq\ 0
\right \}
\ =\ O \left( n^{-\gamma} \right).
\label{eqn-proof-Lemma6.5-c}	
\end{equation}
Notice the reciprocity of (\ref{eqn-proof-Lemma6.5-c}) to (\ref{eqn-proof-Lemma7.4-b-2-aa}).  
%Next, we present the following lemma whose proof is deferred to the end.
%\begin{lemma}  
%\label{proof-Thm5.1-forward-lem4}  For any $\,\gamma > 2,$ when the parameter $\,a\,$ in (\ref{eqn-sparse-exploration}) is selected to be sufficiently close to 1, it holds for every $\,i = 0, \ldots, M-1,$ and any $\,\epsilon' > 0\,$ that the first time $\,T\,$ after which the ML estimate $\,\hat{i}\,$ always equals the true hypothesis $\,i,$ satisfies
%\begin{equation}
%\lim\limits_{n \rightarrow \infty}
%\mathbb{P}_i \left \{
%T \ \geq\ \epsilon' n
%\right \}  \ =\ O \left( n^{- \gamma} \right).
%\label{eqn-proof-Thm5.1-forward-lem4}
%\end{equation}
%\end{lemma}
Applying Lemma \ref{proof-Thm5.1-forward-lem4} with a sufficiently small $\,\epsilon'\,$ yields that to prove (\ref{eqn-proof-Lemma6.5-c}), it suffices to prove that
\begin{equation}
\mathbb{P}_i 
\left \{
\sum\limits_{k=\lfloor \epsilon' n \rfloor}^n
\left (
D \left( p_i^{U_k} \left( \cdot \vert Y_{k-1} \right) \|
				  p_j^{U_k} \left( \cdot \vert Y_{k-1} \right) \right)		 
- 
c \left( U_k \right)
d_i^*
+ \frac{\tilde{\epsilon}}{2}
\right )
\leq 0;\ \tilde{T}  < \epsilon' n
\right \}	\hspace{1.2in} \nonumber	
\end{equation}
% CHECK THIS STEP WHEN WE CHANGE THE MEASURE
\begin{equation}
\hspace{0.2in}   \leq\ 
\tilde{\mathbb{P}}_i 
\left \{
\sum\limits_{k=\lfloor \epsilon' n \rfloor}^n
\left (
D \left( p_i^{U_k} \left( \cdot \vert Y_{k-1} \right) \|
				  p_j^{U_k} \left( \cdot \vert Y_{k-1} \right) \right)		 
- 
c \left( U_k \right)
d_i^*
+ \frac{\tilde{\epsilon}}{2}
\right )
\leq 0
\right \}
= O \left( n^{-\gamma} \right),
\label{eqn-proof-Lemma6.5-d}
\end{equation}
where $\,\tilde{\mathbb{P}}_i\,$ is another joint probability measure of $\,\left( Y_{k-1}, U_k \right),\ k = \lfloor \epsilon' n \rfloor, \ldots, n,$ with the same marginal distribution of $\,Y_{\lfloor \epsilon' n \rfloor-1}$ as in $\,\mathbb{P}_i,$ but with each of the $\,U_k\,$ being conditionally independent of $\,U^{k-1}, Y^{k-2}\,$ conditioned on $\,Y_{k-1} = \tilde{y},$ and being conditionally distributed according to $\,q_i^* \left( \cdot \vert \tilde{y} \right),$ where $\,q_i^*\,$ is the argument maximizer in (\ref{eqn-opt-exploit-control-cost}).  Consequently, under $\,\tilde{\mathbb{P}}_i,\ Y_{k},\ k = \lfloor \epsilon' n \rfloor - 1, \ldots, n-1,$ is a stationary Markov chain with transition probabilities being\\
$
p_i^{q^*_i} \left( y \vert \tilde{y} \right) \ =\ 
\sum\limits_{u \in \mathcal{U}} 
q_i^* \left( u \vert \tilde{y} \right)  p_i^{u} \left( y \vert \tilde{y} \right).$

For any $\,\tilde{y} \in \mathcal{Y},$ let 
\begin{eqnarray}
a\left( \tilde{y} \right)  &=& \sum\limits_{u} q_i^* \left( u \vert \tilde{y} \right) 
				D \left( p_i^{u} \left( \cdot \vert \tilde{y} \right) \|
					  p_j^{u} \left( \cdot \vert \tilde{y} \right) \right)	\\
b\left( \tilde{y} \right)  &=& \sum\limits_{u} q_i^* \left( u \vert \tilde{y} \right) 
				c \left( u \right).	  
\end{eqnarray}
Continuing upper bounding the probability on the right-side of (\ref{eqn-proof-Lemma6.5-d}), we get that
\begin{equation}
\tilde{\mathbb{P}}_i 
\left \{
\sum\limits_{k=\lfloor \epsilon' n \rfloor}^n
\left (
D \left( p_i^{U_k} \left( \cdot \vert Y_{k-1} \right) \|
				  p_j^{U_k} \left( \cdot \vert Y_{k-1} \right) \right)		 
- 
c \left( U_k \right)
d_i^*
+ \frac{\tilde{\epsilon}}{2}
\right )
\leq 0
\right \}
	\hspace{1.5in} \nonumber	
\end{equation}
\begin{equation}
\hspace{0.2in} =\ 
\tilde{\mathbb{P}}_i 
\left \{
\begin{array}{cc}
\sum\limits_{k=\lfloor \epsilon' n \rfloor}^n
\left (
D \left( p_i^{U_k} \left( \cdot \vert Y_{k-1} \right) \|
				  p_j^{U_k} \left( \cdot \vert Y_{k-1} \right) \right)		 
- 
a \left( Y_{k-1} \right)
+ \frac{\tilde{\epsilon}}{6}
\right )  \\
\sum\limits_{k=\lfloor \epsilon' n \rfloor}^n
\left (
b \left( Y_{k-1} \right)
d_i^*
- 
c \left( U_k \right)
d_i^*
+ \frac{\tilde{\epsilon}}{6}
\right )  \\
\sum\limits_{k=\lfloor \epsilon' n \rfloor}^n
\left (
a \left( Y_{k-1} \right)
- 
b \left( Y_{k-1} \right)
d_i^*
+ \frac{\tilde{\epsilon}}{6}
\right )  \\
\end{array}
\ \leq\  0
\right \}.
\label{eqn-proof-Lemma6.5-e}
\end{equation}
That the two probabilities 
\begin {eqnarray}
\tilde{\mathbb{P}}_i 
\left \{
\sum\limits_{k=\lfloor \epsilon' n \rfloor}^n
\left (
D \left( p_i^{U_k} \left( \cdot \vert Y_{k-1} \right) \|
				  p_j^{U_k} \left( \cdot \vert Y_{k-1} \right) \right)		 
- 
a \left( Y_{k-1} \right)
+ \frac{\tilde{\epsilon}}{6}
\right )
\leq 0
\right \},\nonumber \\ 
\tilde{\mathbb{P}}_i 
\left \{
\sum\limits_{k=\lfloor \epsilon' n \rfloor}^n
\left (
b \left( Y_{k-1} \right)
d_i^*
- 
c \left( U_k \right)
d_i^*
+ \frac{\tilde{\epsilon}}{6}
\right )
\leq 0
\right \} \hspace{0.8in}
\nonumber
\end{eqnarray}
go to zero exponentially in $\,n\,$ follow from the same Chernoff bounding argument as in the one in (\ref{eqn-proof-Lemma6.5-b-2}) upon using Assumption (\ref{eqn-assump-pos-transition}) and finiteness of $\,\mathcal{Y},$ and upon noting that 
\begin{eqnarray}
\tilde{\mathbb{E}_i} 
\left [ 
D \left( p_i^{U_k} \left( \cdot \vert Y_{k-1} \right) \|
				  p_j^{U_k} \left( \cdot \vert Y_{k-1} \right) \right)		 
\bigg \vert Y_{k-1} \right ] &=& a \left( Y_{k-1} \right),	\nonumber \\
\tilde{\mathbb{E}_i} \left [ c \left( U_k \right)
 \vert Y_{k-1} \right ] &=& b \left( Y_{k-1} \right).	\nonumber
\end{eqnarray}
To finish the proof of (\ref{eqn-proof-Lemma6.5-c}), we now prove that 
\begin{equation}
\tilde{\mathbb{P}}_i 
\left \{
\sum\limits_{k=\lfloor \epsilon' n \rfloor}^n
\left (
a \left( Y_{k-1} \right) - b \left( Y_{k-1} \right) d_i^*
+ \frac{\tilde{\epsilon}}{6}
\right )
\leq 0
\right \} \ \leq\ O \left( n^{-\gamma} \right).
\label{eqn-proof-Lemma6.5-e-2}
\end{equation}
To this end, we employ the convex analytic method again.  First, note from (\ref{eqn-target-opt-coeffs}) and Assumption (\ref{eqn-assump-pos-transition}) that 
\begin{equation}
d_i^* \ =\
\frac{
\sum\limits_{\tilde{y} \in \mathcal{Y}}
\mu_i^{q_i^*} \left( \tilde{y} \right) a\left( \tilde{y} \right)
}
{
\sum\limits_{\tilde{y} \in \mathcal{Y}}
\mu_i^{q_i^*} \left( \tilde{y} \right) b\left( \tilde{y} \right)
},
\label{eqn-proof-Lemma6.5-f}
\end{equation}
where $\,\mu_i^{q_i^*}\,$ is the unique distribution satisfying for every $\,y \in \mathcal{Y}\,$ that
\begin{equation}
\mu_i^{q_i^*} \left( y \right)  
= \sum\limits_{\tilde{y}\in \mathcal{Y}}
p_i^{q^*_i} \left( y \vert \tilde{y} \right)~\mu_i^{q_i^*} \left( \tilde{y} \right).
\label{eqn-proof-Lemma6.5-f-2}
\end{equation}
Next, since $\, b \left( \tilde{y} \right ) \geq \underline{c},$ we get that 
the quantity $\frac{
\sum\limits_{\tilde{y} \in \mathcal{Y}}
t \left( \tilde{y} \right) a\left( \tilde{y} \right)
}
{
\sum\limits_{\tilde{y} \in \mathcal{Y}}
t \left( \tilde{y} \right) b\left( \tilde{y} \right)
}\,$ is a continuous function of the distribution $\,t$.  By this continuity and the fact that $\,\mu_i^{q_i^*}\,$ is the unique solution to (\ref{eqn-proof-Lemma6.5-f-2}), it holds that there exists a sufficiently small $\,\epsilon\,$ such that
\begin{equation}
d_i^* - \frac{\tilde{\epsilon}}{6 \overline{c}}  \ \ <\ \ \ 
\min\limits_{t:\ 
\bigg \vert t \left( y \right)  
- \sum\limits_{\tilde{y}\in \mathcal{Y}}
p_i^{q^*_i} \left( y \vert \tilde{y} \right)~t \left( \tilde{y} \right) 
\bigg \vert \ \leq\  \epsilon,\ \forall y \in \mathcal{Y}
} \ \ \  \frac{
\sum\limits_{\tilde{y} \in \mathcal{Y}}
t \left( \tilde{y} \right) a\left( \tilde{y} \right)
}
{
\sum\limits_{\tilde{y} \in \mathcal{Y}}
t \left( \tilde{y} \right) b\left( \tilde{y} \right)
}.
\label{eqn-proof-Lemma6.5-g}
\end{equation}
Now we look at the empirical distribution
\begin{equation}
	T \left( \tilde{y} \right) \ =\
	\frac{1}{n - \lfloor \epsilon' n \rfloor + 1} 
	\sum\limits_{k={\lfloor \epsilon' n \rfloor}}^n  \mathbb{I} \left \{ Y_{k-1} = \tilde{y} \right \}.
	\nonumber	
\end{equation}
Then, if follows from (\ref{eqn-proof-Lemma6.5-g}) that
\begin{equation}
\tilde{\mathbb{P}}_i 
\left \{
\sum\limits_{k=\lfloor \epsilon' n \rfloor}^n
\left (
a \left( Y_{k-1} \right) - b \left( Y_{k-1} \right) d_i^*
+ \frac{\tilde{\epsilon}}{6}
\right )
\leq 0
\right \}	\hspace{2.0in}
\nonumber
\end{equation}
\begin{eqnarray} 
&\leq&
\tilde{\mathbb{P}}_i 
\left \{
\frac{\sum\limits_{\tilde{y} \in \mathcal{Y}} T \left( \tilde{y} \right) a \left( \tilde{y} \right)}
{\sum\limits_{\tilde{y} \in \mathcal{Y}} T \left( \tilde{y} \right) b \left( \tilde{y} \right)}
\leq d_i^* - \frac{\tilde{\epsilon}}{6 \overline{c}}
\right \} 	\nonumber \\
&\leq&
\tilde{\mathbb{P}}_i 
\left \{
\bigg \vert T \left( y \right)  
- \sum\limits_{\tilde{y}\in \mathcal{Y}}
p_i^{q^*_i} \left( y \vert \tilde{y} \right)~T \left( \tilde{y} \right) 
\bigg \vert > \epsilon,\ \exists y \in \mathcal{Y}
\right \}.
\label{eqn-proof-Lemma6.5-h}
\end{eqnarray}
The same Chernoff bounding argument applied to 
\begin{equation}
\sum\limits_{k={\lfloor \epsilon' n \rfloor}}^n  
\left(
\mathbb{I} \left \{ Y_{k} = y \right \}
- \tilde{\mathbb{E}}_i \left [
\mathbb{I} \left \{ Y_{k} = y \right \}
\big \vert Y^{k-1}
\right ]
\right)
\hspace{2.0in}
\nonumber 
\end{equation}
\begin{equation}
\ \ \ = 
\sum\limits_{k={\lfloor \epsilon' n \rfloor}}^n  
\left(
\mathbb{I} \left \{ Y_{k} = y \right \}
- \sum\limits_{\tilde{y} \in \mathcal{Y}}
p_i^{q^*_i} \left( y \vert \tilde{y} \right) 
~\mathbb{I} \left \{ Y_{k-1} = \tilde{y} \right \}
\right)
\nonumber
\end{equation}
yields that the probability on the right-side of (\ref{eqn-proof-Lemma6.5-h}) vanishes exponentially in $\,n,$ thereby proving (\ref{eqn-proof-Lemma6.5-e-2}).
%CHECK THIS STEP

Lastly, it is only left to prove Lemma \ref{proof-Thm5.1-forward-lem4}.
%which also follows from the same Chernoff bounding argument as in the proof of Theorem 3 in Appendix B.II of \cite*{niti-atia-veer-ieeetac-2013}.  
We note that for each $\,j \neq i,\ -1 < t < 0,\ \ell \geq \lfloor \epsilon' n \rfloor -1,$
\begin{eqnarray}
\mathbb{P}_i \left \{
\sum\limits_{k=1}^{\ell}
	\log{\left( \frac{p_i^{U_k} \left( Y_k \vert Y_{k-1} \right)}
	{p_j^{U_k} \left( Y_k \vert Y_{k-1} \right) } \right)}
\ \leq\ 0
\right \}
&\leq&
\mathbb{E}_i \left [ 
e^{ 
t
\left (
\sum\limits_{k=1}^{\ell}
	\log{\left( \frac{p_i^{U_k} \left( Y_k \vert Y_{k-1} \right)}
	{p_j^{U_k} \left( Y_k \vert Y_{k-1} \right) } \right)}
\right )
}
\right ].
\label{eqn-Pf-Lemma6.6-a} 
\end{eqnarray}
Note that by the convexity of the moment generating function, 
it holds for every $\,k = 1, \ldots, \ell\,$ that
\begin{equation}
\mathbb{E}_i \left [ 
e^{ 
t
\left (
	\log{\left( \frac{p_i^{U_k} \left( Y_k \vert Y_{k-1} \right)}
	{p_j^{U_k} \left( Y_k \vert Y_{k-1} \right) } \right)}
\right)
}
\bigg \vert Y^{k-1}, U^{k-1}
\right ] \ \leq\  1 \ \mbox{a.s.}.	\nonumber
\end{equation}
Also by (\ref{eqn-assump-pos-transition}) and finiteness of $\,\mathcal{U}, \mathcal{Y},$ we get that at all those times $\,k\,$ in (\ref{eqn-sparse-exploration}) wherein the control $\,U_k\,$ is picked to be uniformly distributed, 
\begin{equation}
\mathbb{E}_i \left [ 
e^{ 
t
\left (
	\log{\left( \frac{p_i^{U_k} \left( Y_k \vert Y_{k-1} \right)}
	{p_j^{U_k} \left( Y_k \vert Y_{k-1} \right) } \right)}
\right)
}
\bigg \vert Y^{k-1}, U^{k-1}
\right ] \ \leq\  e^{-b'} < 1 \ \mbox{a.s.},	\nonumber
\end{equation} 
for some $\,b' > 0.$  Hence, we get from (\ref{eqn-Pf-Lemma6.6-a}) that
\begin{equation}
\mathbb{P}_i \left \{ \tilde{T} \geq \epsilon' n \right \} \leq
\sum\limits_{\ell = \lfloor \epsilon' n \rfloor -1}^{\infty}
\mathbb{P}_i \left \{
p_i \left( Y^{\ell}, U^{\ell} \right)
\leq p_j \left( Y^{\ell}, U^{\ell} \right)
\right \} \leq 
\sum\limits_{\ell = \lfloor \epsilon' n \rfloor -1}^{\infty} 
O \left( e^{-b'  \frac{\log{\ell}}{\log{a}} } \right) = 
O \left( n^{-\gamma} \right),
\nonumber
\end{equation}
for a suitably chosen $\,a\,$ in (\ref{eqn-sparse-exploration}) close enough to 1, thereby completing the proof of the lemma.

\vspace{0.2in}
\noindent
{\em Proof That the Test with the Stopping Rule (\ref{eqn-stopping-hardriskconstraints}) Meets Distinctly Predefined Risk Constraints}
\vspace{0.1in}

%The proof is similar to the proof of \cite*{niti-atia-veer-ieeetac-2013}[Theorem 4].  
We first prove (\ref{eqn-hard-risk-constraints}).  For any $\,j = 0, \ldots, M-1,$ consider the event $\,A_{n,j} = \left \{ \tdblue{N^*} = n,\ \delta = j \right \}$.  It now follows from the stopping rule (\ref{eqn-stopping-hardriskconstraints}) that on $\,A_{n,j},$ we get that
\tdblue{
\begin{eqnarray}
\frac{p_j \left( y^n, u^n \right)}{p_i \left( y^n, u^n \right)}
&\geq&
\left( \frac{p_j \left( y^n, u^n \right)}
{\max\limits_{i \neq j} p_i \left( y^n, u^n \right)} \right)	\nonumber \\
&\geq& 
\frac{1}{\bar{R}_j}.
\end{eqnarray}
}
Consequently, it holds that
\tdblue{
\begin{eqnarray}
\mathbb{P}_i \left \{ \delta = j \right \}
&=&	\sum\limits_{n=1}^{\infty} \mathbb{P}_i \left \{ A_{n, j} \right \}		\nonumber \\
&\leq& \bar{R}_j
	\sum\limits_{n=1}^{\infty} \mathbb{P}_j \left \{ A_{n, j} \right \}		\nonumber \\
&\leq& \bar{R}_j,	\label{eqn-PfThm4.2-bd-error-ij} 
\end{eqnarray}
}
and, hence, 
\tdblue{
\begin{equation}
\tdblue{R_j^*} \ =\ \max\limits_{i \neq j} \mathbb{P}_i \left \{ \delta = j \right \}
\ \leq\ \bar{R}_j. 	\nonumber
\end{equation}
}

It is now left to prove that this test in Section \ref{ss:TestControlCost} but with the stopping rule (\ref{eqn-stopping-hardriskconstraints}) achieves (\ref{eqn-achievable-perf-cost1}) and (\ref{eqn-achievable-perf-cost2}).  To this end, we consider yet another fictitious test following the same control policy as in the previous test, but with the stopping rule (\ref{eqn-stopping-hardriskconstraints}) being replaced by the following rule with a single threshold
\begin{equation}
N' \ \triangleq\ \min\limits_{n \geq 1}~
	\frac{p_{\hat{i}_n} \left(Y^n, U^n \right)}
	{\max\limits_{j \neq \hat{i}_n}~p_j \left( Y^n, U^n \right)}
	~>~
	 \tdblue{T \triangleq \frac{ 1}{\min\limits_{i}  \bar{R}_{i}}. }
	\label{eqn-pf-Thm4.2-fict-test}
\end{equation}
It now follows from (\ref{eqn-stopping-hardriskconstraints}) and (\ref{eqn-pf-Thm4.2-fict-test}) that with probability 1,
\begin{equation}
	\tdblue{N^*} \ \leq\ N'.	\nonumber
\end{equation}
Since we stick with the same control policy in this new test, we also get that with probability 1,
\begin{equation}
	\sum\limits_{k=1}^{\tdblue{N^*}} c\left( U_k \right) \ \leq\ 
		\sum\limits_{k=1}^{N'} c\left( U_k \right).
		\label{eqn-pf-Thm4.2-cost-dom}
\end{equation}
On the other hand, it follows from (\ref{eqn-PfThm4.2-bd-error-ij}) and the assumption in Theorem \ref{thm-optasymperf-nocost-hardrisk} that
\begin{equation}
\tdblue{P_{ {\scriptsize \mathrm{max}}}^*}
\ \leq\ (M-1) \max\limits_{i \neq j} \mathbb{P}_i \left \{ \delta = j \right \}
\ \leq\ \tdblue{(M-1) \max\limits_{i} \bar{R}_i}
\ \leq\ K' \frac{1}{T}.
\label{eqn-pf-Thm4.2-Perror-dom}
\end{equation}

As $\,\max\limits_{i = 0, \ldots, M-1} \bar{R}_i \rightarrow 0,$ the single threshold $\,T\,$ on the right-side of (\ref{eqn-pf-Thm4.2-fict-test}) will go to infinity.  We now get from (\ref{eqn-proof-Thm5.1-forward-b}) in the proof of Theorem \ref{thm-optasymperf-cost} that
\begin{equation}
\lim\limits_{P_{ {\scriptsize \mathrm{max}}} \rightarrow 0}~ 
\frac{\mathbb{E}_i \left [ \sum\limits_{k=1}^{N'} c\left( U_k \right) \right ]}{\log{T}} 
\ \leq\ 
\frac{1}{
\max\limits_{q \left( \cdot \vert \cdot \right)}
~\frac{\min\limits_{j \neq i}
~\sum\limits_{\tilde{y} \in \mathcal{Y}, u \in \mathcal{U}}
\mu_i^q \left( \tilde{y} \right) q \left( u \vert \tilde{y} \right)
D \left( p_i^u \left( \cdot \| \tilde{y} \right) \| p_j^u \left( \cdot \| \tilde{y} \right) \right)}
{\sum\limits_{\tilde{y} \in \mathcal{Y}, u \in \mathcal{U}} \mu_i^q \left( \tilde{y} \right) q \left( u \vert \tilde{y} \right) c(u)}
}.
\label{eqn-pf-Thm4.2-N'-asympperf}
\end{equation}

Lastly, we get from (\ref{eqn-pf-Thm4.2-cost-dom}), (\ref{eqn-pf-Thm4.2-Perror-dom}) and (\ref{eqn-pf-Thm4.2-N'-asympperf}) that 
\begin{eqnarray}
\lim\limits_{P_{ {\scriptsize \mathrm{max}}} \rightarrow 0}~ 
\frac{\mathbb{E}_i \left [ \sum\limits_{k=1}^{\tdblue{N^*}} c \left( U_k \right) \right ]}
{-\log{\left( 
\tdblue{P_{ {\scriptsize \mathrm{max}}}^*}
\right)}} 
&\leq&
\lim\limits_{P_{ {\scriptsize \mathrm{max}}} \rightarrow 0}~ 
\frac{\mathbb{E}_i \left [ \sum\limits_{k=1}^{N'} c \left( U_k \right) \right ]}{\log{T}} 	\nonumber \\
&\leq&
\frac{1}{
\max\limits_{q \left( \cdot \vert \cdot \right)}
~\frac{\min\limits_{j \neq i}
~\sum\limits_{\tilde{y} \in \mathcal{Y}, u \in \mathcal{U}}
\mu_i^q \left( \tilde{y} \right) q \left( u \vert \tilde{y} \right)
D \left( p_i^u \left( \cdot \| \tilde{y} \right) \| p_j^u \left( \cdot \| \tilde{y} \right) \right)}
{\sum\limits_{\tilde{y} \in \mathcal{Y}, u \in \mathcal{U}} \mu_i^q \left( \tilde{y} \right) q \left( u \vert \tilde{y} \right) c(u)}
}.
\nonumber
\end{eqnarray}

%%% Venu's Comments
%%% 1.  In Proposition 2.1 for the uncontrolled observations, we can also do it for stationary Markov sources if we want to keep it consistent with the rest of the paper (see Dragalin-Tatarkovsky-Veeravalli)
%%% 2.  In the discussion section refers to the work of Kiefer and Sacks and say that they have a more practical causal control policy which try to minimize the switching of implementation of instantaneous control value.

%%% Checks
%%% 1. That Proposition 1 works when some divergence is infinite, but all the minimum divergences are finite
%%% 2. Check that in our MSPRT the difference that there is a maximum in the denominator instead of the sum as in Baum-Veeravalli would not change its asymptotic performance and its ability to meet the hard risk constraints
%%% 3. Check that in the models with control cost, the cost function should only affect the control policy and not the stopping rule.
%%% 4. Find a connection between the scheme of Kiefer and Sacks and the cost-biased scheme in the first paper of Becker and Kumar.


\begin{thebibliography}{99}

\setlength{\parskip}{-1.5ex plus0ex minus1.1ex}

\bibitem[Anscombe, 1963]{ansc-jasta-1963}
Anscombe, F. J. (1963). Sequential Medical Trials, {\em Journal of the American Statistical Association} 58: 365--383.

\bibitem[Banerjee and Veeravalli, 2012]{bane-veer-sqa-2012}
Banerjee, T., and Veeravalli, V. V. (2012).  Data-Efficient Quickest Change Detection
with On-Off Observation Control, {\em Sequential Analysis,} 31: 40-77.

\bibitem[Baum and Veeravalli, 1994]{baum-veer-ieeetit-1994}
Baum, C. W., and Veeravalli, V. V. (1994).  A Sequential Procedure for Multiple Hypothesis Testing, {\em IEEE Transactions on Information Theory} 40: 1994--2007.

\bibitem[Bessler, 1960]{bess-tech-repo-1960}
Bessler, S. A. (1960).  Theory and Applications of the Sequential Design of Experiments, $k$-Actions and Infinitely Many Experiments:~Part 1-Theory, Technical report no. 55, Department of Statistics, Stanford University.

\bibitem[Borkar, 1991]{bork-cmc-book-1991}
Borkar, V. S. (1991).  {\em Topics in Controlled Markov Chains,}  Pitman Research Notes in Math. No. 240, Essex, UK: Longman Scientific and Technical.

\bibitem[Burna\v{s}ev, 1976]{burn-prob-per-inf-1976}
Burna\v{s}ev, M. V. (1976).  Data Transmission over a Discrete Channel with Feedback. Random Transmission Time, {\em Problemy Peredachi Informatsii} 12: 10--30.

\bibitem[Chernoff, 1959]{cher-amstat-1959}
Chernoff, H. (1959).  Sequential Design of Experiments, {\em Annals of Mathematical Statistics} 30: 755--770.

\bibitem[Como {\em et al.}, 2009]{como-yuks-tati-ieeetit-2009}
Como, G., Y\"uksel S., and Tatikonda, S. (2009).  The Error Exponent of Variable-Length Codes over Markov Channels with Feedback, {\em IEEE Transactions on Information Theory} 55: 2139--2160.

\bibitem[Dragalin {\em et al.}, 1999]{drag-et-al-ieeetit-1999}
Dragalin, V. P., Tartakovsky, A. G., and Veeravalli, V. V. (1999).  Multihypothesis Sequential Probability Ratio Tests--Part \mbox{I}: Asymptotic Optimality, {\em IEEE Transactions on Information Theory} 45: 2448--2461.

\bibitem[Fuemmeler and Veeravalli, 2008]{fuem-veer-ieeetsp-2008}
Fuemmeler, J. A., and Veeravalli, V. V. (2008).  Smart Sleeping Policies for Energy Efficient Tracking in Sensor Networks, {\em IEEE Transactions on Signal Processing} 56: 2091--2101.

\bibitem[Hero {\em et al.}, 2008]{hero-etal-sens-book-2008}  Hero III, A. O., Casta{\~n}{\'o}n, D. A., Cochran D., and Kastella K. (Eds.) (2008).  {\em Foundations and Applications of Sensor Management,}  Springer Series on Signals and Communication Technology, Springer.

\bibitem[Kiefer and Sacks, 1963]{kief-sack-amstat-1963} Kiefer, J., and Sacks, J. (1963).  Asymptotically Optimum Sequential Inference and Design, {\em Annals of Mathematical Statistics} 34: 705--750.

\bibitem[Kumar and Becker, 1982]{kuma-beck-ieeetac-1982} Kumar, P. R., and Becker, A. (1982).  A New Family of Optimal Adaptive Controllers for Markov Chains, {\em IEEE Transactions on Automatic Control} 27: 137--146. 

\bibitem[Kumar and Varaiya, 1986]{kuma-vara-ss-book-1986}
Kumar, P. R., and Varaiya, P. (1986).  {\em Stochastic Systems: Estimation, Identification, and Adaptive Control,} New Jersey: Prentice-Hall, Inc.

\bibitem[Lo\`{e}ve, 1977]{loev-prob-book-1977}
Lo\`{e}ve, M. (1977). {\em Probability Theory II,} 4th ed., New York: Springer.

\bibitem[Naghshvar and Javidi, 2013]{nagh-javi-astat-2013}
Naghshvar, M. and Javidi, T. (2013).  Active Sequential Hypothesis Testing, {\em Annals of Statistics,} \tdblue{41: 2703-2738.}

\bibitem[Nakibo\u{g}lu and Gallager, 2008]{naki-gall-ieeetit-2008}
Nakibo\u{g}lu, B., and Gallager, R. G. (2008).  Error Exponents for Variable-Length Block Codes with Feedback and Cost Constraints,  
{\em IEEE Transactions on Information Theory} 54: 945--963.
    
\bibitem[Nitinawarat {\em et al.}, 2013]{niti-atia-veer-ieeetac-2013}
Nitinawarat, S., Atia, G., and Veeravalli, V. V. (2013).  Controlled Sensing for Multihypothesis Testing, {\em IEEE Transactions on Automatic Control,} \tdblue{58: 2451-2464}.

\bibitem[Ooi and Wornell, 1998]{ooi-worn-ieeetit-1998}
Ooi, J. M, and Wornell, G. W. (1998).  Fast Iterative Coding Techniques for Feedback Channels,  {\em IEEE Transactions on Information Theory} 44: 2960--2976.


\bibitem[Stone, 1989]{ston-search-book-1989}
Stone, L. D. (1989).  {\em Theory of Optimal Search,}  
Maryland: INFORMS.

\bibitem[Veeravalli and Baum, 1995]{veer-baum-ieeetit-1995}
Veervalli, V. V., and Baum, C. W. (1995).  Asymptotic Efficiency of a Sequential Multihypothesis Testing, {\em IEEE Transactions on Information Theory} 41: 1994--1997.

\bibitem[Yamamoto and Itoh, 1979]{yama-itoh-ieeetit-1979}
Yamamoto, H., and Itoh, K. (1979).  Asymptotic Performance of a Modified Schalkwijk-Barron Scheme for Channels with Noiseless Feedback, {\em IEEE Transactions on Information Theory} 25: 729--733.   

%\bibitem[Lai(1995)]{LaiJRSS95}
%Lai, T. L. (1995). Sequential Changepoint Detection in Quality Control and
%Dynamical Systems, {\em Journal of Royal Statistal Society, Series B} 57: 613--658.

%\bibitem[Siegmund(1985)]{Siegmundbook85}
%Siegmund, D. (1985). {\em Sequential Analysis: Tests and Confidence Intervals}, New
%York: Springer-Verlag.

%\bibitem[Shiryaev(1963)]{ShiryaevTPA63}
%Shiryaev, A. N. (1963). On Optimum Methods in Quickest Detection Problems,
%{\em Theory of Probability and Its Applications} 8: 22--46.

%\bibitem[Shiryaev(1978)]{ShiryaevBook78}
%Shiryaev, A. N. (1978). {\em Optimal Stopping Rules}, New York: Springer-Verlag.

%\bibitem[Tartakovsky(2005)]{TartakovskyIEEECDC05}
%Tartakovsky, A. G. (2005). Asymptotic Performance of a Multichart CUSUM Test Under
%False Alarm Probability Constraint, \emph{Proceedings of 44th IEEE Conference on Decision
%and Control and European Control Conference (CDC-ECC'05)}, December 12-15,
%2005, pp.\ 320--325, Seville, Spain, Omnipress CD-ROM, ISBN 0-7803-9568-9.

%\bibitem[Tartakovsky et al.(2008)]{TartakovskyPollakPolIWAP08}
%Tartakovsky, A. G., Pollak, M. and Polunchenko, A. (2008). Asymptotic
%Exponentiality of First Exit Times for Recurrent Markov Processes and Applications
%to Changepoint Detection, {\em International Workhop on Applied Probability}, July
%7-10, 2008, Compie\'gne, France.

%\bibitem[Tartakovsky and Veeravalli(2004)]{TarVeerASM2004}
%Tartakovsky, A. G. and Veeravalli, V. V. (2004). Change-Point Detection in
%Multichannel and Distributed Systems with Applications, {\em Applications of
%Sequential Methodologies}, N.\ Mukhopadhyay, S.\ Datta and S.\ Chattopadhyay, eds.,
%pp.\ 339--370, New York: Marcel Dekker.

%\bibitem[Tartakovsky and Veeravalli(2005)]{TartakovskyVeerTVP05}
%Tartakovsky, A. G. and Veeravalli, V. V. (2005). General Asymptotic Bayesian Theory
%of Quickest Change Detection, {\em Theory of Probability and Its Applications} 49:
%458--497.

\end{thebibliography}
\end{document}